\documentclass[12pt,a4paper]{tac}
\usepackage{amssymb, amsmath, stmaryrd, mathrsfs}
\usepackage[latin1]{inputenc}
\usepackage[arrow, matrix, tips, curve]{xy}
\usepackage{hyperref}
\SelectTips{cm}{10}
\usepackage{enumitem}
\setlist[enumerate]{noitemsep}
\setlist[enumerate,1]{label=(\alph*), ref=(\alph*)}
\setlist[enumerate,2]{label=(\roman*),
ref=\theenumi(\roman*)}

\author{John Bourke and Nick Gurski}

\thanks{The first author acknowledges the support of the Grant agency
  of the Czech Republic, grant number P201/12/G028.  The second author acknowledges the support of the EPSRC in the form of an Early Career Fellowship.  Both authors thank the University of Sheffield for support via the Mathematical Sciences Research Centre.}

\address{Department of Mathematics and Statistics, Masaryk University\\ Kotl\'a\v rsk\'a 2, Brno 60000, Czech Republic\\[5pt]
School of Mathematics and Statistics, University of Sheffield\\ Sheffield, United Kingdom\\
}

\title{A cocategorical obstruction to tensor products of Gray-categories}

\copyrightyear{2015}


\eaddress{bourkej@math.muni.cz\CR nick.gurski@sheffield.ac.uk}

\makeatletter

\renewcommand{\paragraph}{\@startsection
{paragraph}%
{3}%
{0mm}%
{-\baselineskip}%
{-0.4em plus 0.2em minus 0.2em}%
{\normalfont\normalsize\bfseries}}%
\makeatother

\makeatletter \@namedef{itemize*}{\itemize\parsep\z@ \parskip\z@}
\@namedef{enditemize*}{\enditemize} \@namedef{enumerate*}{\enumerate\parsep\z@
\parskip\z@} \@namedef{endenumerate*}{\endenumerate}
\makeatother
\DeclareMathAlphabet      {\mathbf}{OT1}{cmr}{b}{n}

\let\pf\proof
\let\epf\endproof

\numberwithin{equation}{section}
\theoremstyle{plain}

\newtheorem{Theorem}{Theorem}[section]

\newtheorem{Corollary}[Theorem]{Corollary}
\newtheorem{Proposition}[Theorem]{Proposition}
\newtheorem{Lemma}[Theorem]{Lemma}

\def\proofof[#1]{\subsubsection*{{\textsc{Proof of #1.}}}}

\theoremstyle{definition}

\newcommand{\vcat}{\ensuremath{\mathcal{V}\textnormal{-Cat}}\xspace}
\newcommand{\twocat}{\ensuremath{\textnormal{2-Cat}}\xspace}
\newcommand{\OCat}{\ensuremath{\omega \textnormal{-Cat}}\xspace}
\newcommand{\Cat}{\ensuremath{\textnormal{Cat}}\xspace}
\newcommand{\Der}{\ensuremath{\textnormal{Der}}\xspace}
\newcommand{\GCat}{\ensuremath{\textnormal{Gray-Cat}}\xspace}
\newcommand{\SCat}{\textnormal{Sesquicat}\xspace}
\newcommand{\Set}{\ensuremath{\textnormal{Set}}\xspace}
\newcommand{\f}[1]{\ensuremath{\mathcal{#1}}\xspace}
\newcommand{\arr}{\ensuremath{\mathbf{S}}\xspace}
\newcommand{\Ob}{\ensuremath{\mathbf{O}}\xspace}
\newcommand{\cd}[2][]{\vcenter{\hbox{\xymatrix#1{#2}}}}
\newcommand{\atwo}{\textbf{2}\xspace}
\newcommand{\athree}{\textbf{3}\xspace}

\begin{document}

\maketitle

\begin{abstract}
It was argued by Crans that it is too much to ask that the category of Gray-categories admit a well behaved monoidal biclosed structure.  We make this precise by establishing undesirable properties that any such monoidal biclosed structure must have.  In particular we show that there does not exist any tensor product making the model category of Gray-categories into a monoidal model category.
\end{abstract}

\section{Introduction}
The category $\twocat$ of small 2-categories admits several monoidal biclosed structures.  One of these, usually called the Gray-tensor product $\otimes_{p}$, has a number of appealing features.
\begin{enumerate}
\item[(1)] Given 2-categories $A$ and $B$ the morphisms of the corresponding 2-category $[A,B]$ are \emph{pseudonatural transformations}: the most important transformations in 2-category theory.
\item[(2)] Each tricategory is equivalent to a \emph{Gray-category} \cite{Gordon1995Coherence}: a category enriched in $(\twocat,\otimes_{p})$.
\item[(3)] The Gray tensor product equips the model category $\twocat$ with the structure of a \emph{monoidal model category} \cite{Lack2002A-quillen}; in particular, it equips the homotopy category of $\twocat$ with the structure of a monoidal closed category.
\end{enumerate}
Gray-categories obtain importance by virtue of (2).  Instead of working in a general tricategory it suffices to work in a Gray-category -- for an example of this see \cite{Day1997Monoidal}.  They are more manageable than general tricategories, and differ primarily from strict 3-categories only in that the \emph{middle four interchange} does not hold on the nose, but rather up to coherent isomorphism.\\
Bearing the above in mind, it is natural to ask whether there exists a tensor product of Gray-categories satisfying some good properties analogous to the above ones.  This topic was investigated by Crans, who in the introduction to \cite{Crans1999A-tensor} claimed that it is too much to ask for a \emph{monoidal biclosed structure} on $\GCat$ that captures weak transformations, and asserts that this is due to the failure of the middle four interchange.\\
The folklore argument against the existence of such a monoidal biclosed structure --  discussed in \cite{nlab2009} and likely originating with James Dolan -- is fairly compelling and goes as follows.  (For simplicity we here confine ourselves to the symmetric monoidal case.)  Given Gray-categories $A$ and $B$ such a structure would involve a Gray-category $[A,B]$ whose objects should be \emph{strict} Gray-functors\begin{footnote}{In fact biclosedness forces the objects to be the strict Gray-functors -- see Proposition~\ref{prop:unit} and Remark~\ref{rk:HigherUnit}.}\end{footnote} and whose 1-cells $\eta:F \to G$ should be \emph{weak transformations}.  Such a transformation should involve at least 1-cell components $\eta_{a}:Fa \to Ga$ and 2-cells
$$\xy
(0,0)*+{Fa}="00"; (15,0)*+{Fb}="10";
(0,-15)*+{Ga}="01";(15,-15)*+{Gb}="11";
{\ar^{F\alpha} "00"; "10"};
{\ar_{\eta_{a}} "00"; "01"};
{\ar^{\eta_{b}} "10"; "11"};
{\ar_{G\alpha} "01"; "11"};
{\ar@{=>}^{\eta_{\alpha}}(11,-5)*+{};(5,-10)*+{}};
\endxy$$
satisfying coherence equations likely involving 3-cells.  Which equations should be imposed?\newline{}
If we view $\eta:F \to G \in [A,B]$ as a Gray-functor $\atwo \to [A,B]$ it should correspond, by closedness, to a Gray-functor $\overline{\eta}:A \to [\atwo,B]$ where $\overline{\eta}(a)=\eta_{a}:Fa \to Ga$ whilst $\overline{\eta}(\alpha)$ should encode the above 2-cell $\eta_{\alpha}$; the requirement that the Gray-functor $\overline{\eta}:A \to [\atwo,B]$ preserve composition should then correspond to asking that the condition
\begin{equation}
\label{eq:comp}
\begin{gathered}
\xy
(0,0)*+{Fa}="00"; (15,0)*+{Fb}="10";(30,0)*+{Fc}="20";
(0,-15)*+{Ga}="01";(15,-15)*+{Gb}="11";(30,-15)*+{Gc}="21";
{\ar^{F\alpha} "00"; "10"};
{\ar^{F\beta} "10"; "20"};
{\ar_{\eta_{a}} "00"; "01"};
{\ar^{\eta_{b}} "10"; "11"};
{\ar^{\eta_{c}} "20"; "21"};
{\ar_{G\alpha} "01"; "11"};
{\ar_{G\beta} "11"; "21"};
{\ar@{=>}^{\eta_{\alpha}}(11,-5)*+{};(5,-10)*+{}};
{\ar@{=>}^{\eta_{\beta}}(26,-5)*+{};(20,-10)*+{}};
\endxy
\hspace{0.5cm}
\xy
(0,-7)*+{=}="00";
\endxy
\hspace{0.5cm}
\xy
(0,0)*+{Fa}="00"; (20,0)*+{Fc}="10";
(0,-15)*+{Ga}="01";(20,-15)*+{Gc}="11";
{\ar^{F(\beta \alpha)} "00"; "10"};
{\ar_{\eta_{a}} "00"; "01"};
{\ar^{\eta_{c}} "10"; "11"};
{\ar_{G(\beta \alpha)} "01"; "11"};
{\ar@{=>}^{\eta_{\beta \alpha}}(12,-5)*+{};(6,-10)*+{}};
\endxy
\end{gathered}
\end{equation}
holds for $\eta$.  Now given $\eta:F \to G$ and $\mu:G \to H$ satisfying \eqref{eq:comp} we expect to define the composite $\mu \circ \eta:F \to H \in [A,B]$ to have component $(\mu \circ \eta)_{\alpha}$ as left below
$$\xy
(0,0)*+{Fa}="00"; (15,0)*+{Fb}="10";
(0,-15)*+{Ga}="01";(15,-15)*+{Gb}="11";
(0,-30)*+{Ha}="02";(15,-30)*+{Hb}="12";
{\ar^{F\alpha} "00"; "10"};
{\ar_{\eta_{a}} "00"; "01"};
{\ar^{\eta_{b}} "10"; "11"};
{\ar_{G\alpha} "01"; "11"};
{\ar@{=>}^{\eta_{\alpha}}(11,-5)*+{};(5,-10)*+{}};
{\ar_{\mu_{a}} "01"; "02"};
{\ar^{\mu_{b}} "11"; "12"};
{\ar_{H\alpha} "02"; "12"};
{\ar@{=>}^{\mu_{\alpha}}(11,-20)*+{};(5,-25)*+{}};
\endxy
\hspace{2cm}
\xy
(0,0)*+{Fa}="00"; (15,0)*+{Fb}="10";(30,0)*+{Fc}="20";
(0,-15)*+{Ga}="01";(15,-15)*+{Gb}="11";(30,-15)*+{Gc}="21";
(0,-30)*+{Ha}="02";(15,-30)*+{Hb}="12";(30,-30)*+{Hc}="22";
{\ar^{F\alpha} "00"; "10"};
{\ar_{\eta_{a}} "00"; "01"};
{\ar^{\eta_{b}} "10"; "11"};
{\ar_{G\alpha} "01"; "11"};
{\ar@{=>}^{\eta_{\alpha}}(11,-5)*+{};(5,-10)*+{}};
{\ar_{\mu_{a}} "01"; "02"};
{\ar^{\mu_{b}} "11"; "12"};
{\ar_{H\alpha} "02"; "12"};
{\ar@{=>}^{\mu_{\alpha}}(11,-20)*+{};(5,-25)*+{}};
{\ar^{F\beta} "10"; "20"};
{\ar^{\eta_{c}} "20"; "21"};
{\ar_{G\beta} "11"; "21"};
{\ar^{\mu_{c}} "21"; "22"};
{\ar_{H\beta} "12"; "22"};
{\ar@{=>}^{\mu_{\beta}}(26,-20)*+{};(21,-25)*+{}};
{\ar@{=>}^{\eta_{\beta}}(26,-5)*+{};(20,-10)*+{}};
(35,-31)*+{.};
\endxy$$
Having done so, one asks whether $\mu \circ \beta$ satisfies \eqref{eq:comp} and this amounts to the assertion that the two ways of composing the four 2-cells above right -- vertical followed by horizontal and horizontal followed by vertical -- coincide.  But in a general Gray-category they do not.  The conclusion is that that we cannot define a category $[A,B]$ of Gray-functors and weak transformations.
\newline{}
The above argument, though convincing, is not quite precise.  Our goal here is to describe heavy and undesirable restrictions that any monoidal biclosed structure on $\GCat$ must satisfy.  This is the content of our main result, Theorem~\ref{thm:TensorGray}.  This result immediately rules out the possibility of the sort of tensor product that one initially hopes for: in which $\atwo \otimes \atwo$ is the non-commuting square
$$
\xy
(0,0)*+{(0,0)}="00"; (15,0)*+{(1,0)}="10";
(0,-15)*+{(0,1)}="01";(15,-15)*+{(1,1)}="11";
{\ar "00"; "10"};
{\ar "00"; "01"};
{\ar "10"; "11"};
{\ar "01"; "11"};
{\ar@{=>}^{}(11,-3)*+{};(5,-8)*+{}};
{\ar@{=>}^{}(5,-12)*+{};(11,-7)*+{}};
\endxy$$
with some kind of equivalence connecting the two paths $(0,0) \rightrightarrows (1,1)$.  In particular we use Theorem~\ref{thm:TensorGray} to rule out the possibility of a monoidal biclosed structure on $\GCat$ capturing weak transformations or yielding a monoidal model structure.\newline{}
%prove that there can exist no monoidal model structure on $\GCat$.
In establishing these results our line of argument is rather different to that outlined above.  In fact our techniques are based upon \cite{Foltz1980Algebraic}, in which it was shown that \Cat admits exactly two monoidal biclosed structures.  The approach there was to observe that biclosed structures on an algebraic category amount to ``double coalgebras" in the same category, and so to enumerate the former it suffices to calculate the latter, of which there are often but few.  On \Cat this amounts to the enumeration of certain kinds of \emph{double cocategories} in \Cat and we give a careful account of this background material from \cite{Foltz1980Algebraic} in Section 2.  Section 3 recalls the various known tensor products on \twocat.  Our main result is Theorem~\ref{thm:TensorGray} of Section 4.  Corollary~\ref{cor:WeakTransformations} applies this result to describe a precise limitation on the kinds of transformations such a biclosed structure can capture, and Corollary~\ref{cor:model} applies it to rule out the existence of a monoidal model structure.  We conclude by discussing a related argument of James Dolan.\\
The authors would like to thank John Baez, James Dolan, Michael Shulman and Ross Street for useful correspondence concerning a preliminary version of this article.

\section{Two monoidal biclosed structures on \Cat}
It was shown in \cite{Foltz1980Algebraic} that there exist precisely two monoidal biclosed structures on \Cat and an understanding of this result forms the starting point of our analysis of the \GCat situation.  The treatment in \cite{Foltz1980Algebraic} leaves certain minor details to the reader.   Since these details will be important in Section 4, we devote the present section to giving a full account -- but emphasise that nothing here is original.\newline{}
\Cat is cartesian closed and this accounts for one of the monoidal biclosed structures.  Given categories $A$ and $B$ the corresponding internal hom $[A,B]$ is the category of functors and natural transformations.  This forms a subcategory of $[A,B]_{f}$, the category of functors and mere transformations: families of arrows of which naturality is not required.  $[A,B]_{f}$ is the internal hom of the so-called funny tensor product $\star$ \begin{footnote}{The terminology ``funny tensor product" follows \cite{Street1996Categorical}.}\end{footnote}, which also has unit the terminal category $1$.  The funny tensor product $A \star B$ is the pushout $$\xy
(0,0)*+{obA \times obB}="00"; (25,0)*+{obA \times B}="10";
(0,-15)*+{A \times obB}="01"; (25,-15)*+{A\star B}="11";
{\ar^{} "00";"10"};
{\ar ^{}"10";"11"};
{\ar _{}"00";"01"};
{\ar _{} "01";"11"};
\endxy$$
where $obA$ and $obB$ are the discrete categories with the same objects as $A$ and $B$ respectively.  Using this formula one sees that the funny tensor product is symmetric monoidal and, furthermore, that the induced maps $$A \star B \to A \times B$$ from the pushout equip the identity functor $1:\Cat \to \Cat$ with the structure of an opmonoidal functor.
\begin{footnote}{Monoidal categories with unit $1$ are sometimes called semicartesian.  In fact the funny and cartesian tensor products are respectively the initial and terminal semicartesian monoidal structures on \Cat.}\end{footnote}
\subsection{Cocontinuous bifunctors with unit.}
The terminal category $1$ is the unit for both monoidal biclosed structures.  Let us recall why this must be so. \newline{}
Firstly we call a category $\f C$ equipped with a bifunctor $\otimes:\f C \times \f C \to \f C$, object $I$ and natural isomorphisms $l_{A}:I \otimes A \to A$ and $r_{A}:A \otimes I \to A$ such that $l_{I}=r_{I}:I \otimes I \to I$ a \emph{bifunctor with unit}.  If $\f C$ is cocomplete and $\otimes:\f C \times \f C \to \f C$ cocontinuous in each variable we will often refer to it as a \emph{cocontinuous bifunctor} and, when additionally equipped with a unit in the above sense, as a \emph{cocontinuous bifunctor with unit.}
%The notation suggests that the unit is determined by the bifunctor, and in the cases of interest this is true.
\begin{Proposition}\label{prop:unit}
Any cocontinuous bifunctor $\otimes$ with unit on $\Cat$ has unit the terminal category $1$.
\end{Proposition}
\pf
For a category $\f C$, let $End(1_{\f C})$ denote the endomorphism monoid of the identity functor $1_{\f C}$.  In such a setting the morphism of monoids $End(1_{\f C}) \to End(I)$ given by evaluation at $I$ admits a section.  This assigns to $f:I \to I$ the endomorphism of $1_{\f C}$ with component
 $$\cd{A \ar[r]^{r_{A}^{-1}} & A \otimes I \ar[r]^{A \otimes f} & A \otimes I \ar[r]^{r_{A}} & A}$$ at $A$.
The identity functor $1:\Cat \to \Cat$ has a single endomorphism: for given $\alpha \in End(1_{\Cat})$ naturality at each functor $1 \to A$ forces $\alpha_{A}$ to be the identity on objects, whilst naturality at functors $\atwo \to A$ forces $\alpha_{A}$ to be the identity on arrows.  All but the initial and terminal categories admit a non-trivial endomorphism -- a constant functor -- so that the unit $I$ must be either $0$ or $1$.  Cocontinuity of $A \otimes -$ forces $A \otimes 0 \cong 0$ and so leaves $1$ as the only possible unit.\epf
\subsection{Remark.}\label{rk:HigherUnit}
The above argument generalises easily to higher dimensions: in particular, each monoidal biclosed structure on $\twocat$ or $\GCat$ must have unit the terminal object.  The argument in either case extends the above one, additionally using naturality in maps out of the free living $2$-cell/$3$-cell as appropriate.\\   This restriction on the unit forces, in each case, the objects of the corresponding internal hom $[A,B]$ to be the strict 2-functors or Gray-functors respectively.  In the case of $\twocat$ this is illustrated by the natural bijections $\twocat(1,[A,B])\cong \twocat(1 \otimes A,B) \cong \twocat(A,B)$.  The $\GCat$ case is identical.
\subsection{Double coalgebras versus cocontinuous bifunctors.}
Locally finitely presentable categories are those of the form $Mod(T)=Lex(T,\Set)$ for $T$ a small category with finite limits.  A standard reference is \cite{Adamek1994Locally}.  Examples include $\Cat, \twocat$ and $\GCat$.\\
Our interest is in monoidal biclosed structures on such categories.  We note that a monoidal structure $\otimes$ on $Mod(T)$ is biclosed just when the tensor product $- \otimes -$ is cocontinuous in each variable: this follows from the well known fact that each cocontinuous functor $Mod(T) \to \f C$ to a cocomplete category $\f C$ has a right adjoint.  In practice we will use cocontinuity of bifunctors rather than biclosedness.\\
Going beyond $Mod(T)=Lex(T,\Set)$ one can consider the category $Mod(T,\f C)=Lex(T,\f C)$ for any category $\f C$ with finite limits, or if $\f C$ has finite colimits the category $Comod(T,\f C)=Rex(T^{op},\f C)$ of $T$-comodels -- if $T$ is the finite limit theory for categories one obtains the categories $Cat(\f C)$ and $Cocat(\f C)$ of internal categories and cocategories in this way.  The restricted Yoneda embedding $y:T^{op} \to Mod(T)$ preserves finite colimits and is in fact \emph{the universal $T$-comodel in a cocomplete category}: for cocomplete $\f C$ restriction along $y$ yields an equivalence of categories $$Ccts(Mod(T),\f C) \simeq Comod(T,\f C)$$
where on the left hand side $Ccts$ denotes the 2-category of cocomplete categories and cocontinous functors.  At a cocomplete trio let $Ccts(\f A,\f B;\f C)$ denote the category of bifunctors
\[
\f A, \f B \to \f C
\]
cocontinuous in each variable.  Evidently we have an isomorphism $$Ccts(\f A,\f B;\f C)\cong Ccts(\f A,Ccts(\f B,\f C))$$ and applying this together with the above equivalence (twice) gives:

\begin{Proposition}\label{prop:equiv}
For $\f C$ cocomplete the canonical functor $$Ccts(Mod(T),Mod(T);\f C) \to Comod(T,Comod(T,\f C))$$ is an equivalence of categories.
\end{Proposition}
\noindent On the right hand side are $T$-comodels internal to $T$-comodels, or \emph{double $T$-comodels}.

\subsection{Double cocategories.}
We are interested in the special case of the above where $Mod(T) \simeq \Cat$.  Then double $T$-comodels are double cocategories and in this setting we will give a more detailed account.  A cocategory $\mathbf A$ in $\f C$ is a category internal to $\f C^{op}$.  Such consists of a diagram in $\f C$:
\begin{equation}
\label{eq:cocat}
\begin{gathered}
\xy
(0,0)*+{A_{1}}="a0"; (30,0)*+{A_{2}}="b0";(60,0)*+{A_{3}}="c0";
{\ar@<1.5ex>^{d} "a0"; "b0"};
{\ar@<0ex>|{i} "b0"; "a0"};
{\ar@<-1.5ex>_{c} "a0"; "b0"};
{\ar@<1.5ex>^{p} "b0"; "c0"};
%{\ar@<1.5ex>|{r} "c0"; "b0"};
{\ar@<0ex>|{m} "b0"; "c0"};
%{\ar@<-1.5ex>|{l} "c0"; "b0"};
{\ar@<-1.5ex>_{q} "b0"; "c0"};
\endxy
\end{gathered}
\end{equation}
satisfying, to begin with, the reflexive co-graph identities $id=1=ic$.  $A_{3}$ is the pushout $A_{2}+_{A_{1}} A_{2}$ with $pc=qd$, and $m$ the composition map: this satisfies $md=pd$ and $mc=qc$.  The map $m$ is required to be co-associative, a fact we make no use of, and co-unital: this last fact asserts that $m(i,1)=1=m(1,i)$ where $(i,1),(1,i):A_{3} \rightrightarrows A_{2}$ are the induced maps from the pushout characterised by the equations $(i,1)p=di$, $(i,1)q=1$, $(1,i)p=1$ and $(1,i)q=ci$.
\newline{}
The universal cocategory \arr in \Cat, corresponding to the restricted Yoneda embedding, is given by:
\begin{equation}\label{eq:arr}
\begin{gathered}
\xy
%(-20,0)*+{\textnormal{(1.1)}};
(0,0)*+{1}="a0"; (30,0)*+{\atwo}="b0";(60,0)*+{\athree}="c0";
{\ar@<1.5ex>^{d} "a0"; "b0"};
{\ar@<0ex>|{i} "b0"; "a0"};
{\ar@<-1.5ex>_{c} "a0"; "b0"};
{\ar@<1.5ex>^{p} "b0"; "c0"};
{\ar@<0ex>|{m} "b0"; "c0"};
{\ar@<-1.5ex>_{q} "b0"; "c0"};
\endxy
\end{gathered}
\end{equation}
and we sometimes refer to it as the \emph{arrow cocategory}.  It is the restriction of the standard cosimplicial object $\Delta \to \Cat$; in particular $\atwo$ and $\athree$ are the free walking arrow and composable pair:
$$\xy
(0,0)*+{\atwo\, =\, };
(5,0)*+{(0}="a";(15,0)*+{1)}="b";
{\ar^{f} "a"; "b"};
\endxy
\hspace{2cm}
\xy
(0,0)*+{\athree\, =\, };
(5,0)*+{(0}="a";(15,0)*+{1}="b";(25,0)*+{2)}="c";
{\ar^{g} "a"; "b"};
{\ar^{h} "b"; "c"};
(32,-1)*+{.};
\endxy
$$
The various maps involved are order-preserving and characterised by the equations $d < c$ and $p < m < q$ under the pointwise ordering.
\newline{}
A double cocategory is, of course, a cocategory internal to cocategories.  Since colimits of cocategories are pointwise, this amounts to a diagram
 \begin{equation}
 \begin{gathered}
 \label{eq:dbl}
 \xy
%(-20,-30)*+{\textnormal{(1.2)}};
(0,0)*+{A_{1,1}}="00"; (25,0)*+{A_{2,1}}="10";(50,0)*+{A_{3,1}}="20";
%%%%%%%%%%%%%%%
(0,-25)*+{A_{1,2}}="01"; (25,-25)*+{A_{2,2}}="11";(50,-25)*+{A_{3,2}}="21";
%%%%%%%%%%%%%%%%
(0,-50)*+{A_{1,3}}="02"; (25,-50)*+{A_{2,3}}="12";(50,-50)*+{A_{3,3}}="22";
%%%%%%%%%%%%%%%%%%%%%%%%%%%%%%%%%%%%%%%%%%%
{\ar@<2ex>^{d^{1}_{h}} "00"; "10"};
{\ar@<0ex>|{i^{1}_{h}} "10"; "00"};
{\ar@<-2ex>_{c^{1}_{h}} "00"; "10"};
{\ar@<2ex>^{p^{1}_{h}} "10"; "20"};
{\ar@<0ex>|{m^{1}_{h}} "10"; "20"};
{\ar@<-2ex>_{q^{1}_{h}} "10"; "20"};
%%%%%%%%%%%%%%%%%%%%%%%%%%%%%%%%%%%%%%%%%%%%%
{\ar@<2ex>^{d^{2}_{h}} "01"; "11"};
{\ar@<0ex>|{i^{2}_{h}} "11"; "01"};
{\ar@<-2ex>_{c^{2}_{h}} "01"; "11"};
{\ar@<2ex>^{p ^{2}_{h}} "11"; "21"};
{\ar@{.>}@<0ex>|{m^{2}_{h}} "11"; "21"};
{\ar@<-2ex>_{q^{2}_{h}} "11"; "21"};
%%%%%%%%%%%%%%%%%%%%%%%%%%%%%%%%%%%%%%%%%%%%%
{\ar@<2ex>^{d^{3}_{h}} "02"; "12"};
{\ar@<0ex>|{i^{3}_{h}} "12"; "02"};
{\ar@<-2ex>_{c^{3}_{h}} "02"; "12"};
{\ar@<2ex>^{p ^{3}_{h}} "12"; "22"};
{\ar@{.>}@<0ex>|{m^{3}_{h}} "12"; "22"};
{\ar@<-2ex>_{q^{3}_{h}} "12"; "22"};
%%%%%%%%%%%%%%%%%%%%%%%%%%%%%%%%%%%%%%%%%%%%%
{\ar@<2ex>^{c^{1}_{v}} "00"; "01"};
{\ar@<0ex>|{i^{1}_{v}} "01"; "00"};
{\ar@<-2ex>_{d^{1}_{v}} "00"; "01"};
{\ar@<2ex>^{q^{1}_{v}} "01"; "02"};
{\ar@<0ex>|{m^{1}_{v}} "01"; "02"};
{\ar@<-2ex>_{p^{1}_{v}} "01"; "02"};
%%%%%%%%%%%%%%%%%%%%%%%%%%%%%%%%%%%%%%%%%%%
{\ar@<2ex>^{c^{2}_{v}} "10"; "11"};
{\ar@<0ex>|{i^{2}_{v}} "11"; "10"};
{\ar@<-2ex>_{d^{2}_{v}} "10"; "11"};
{\ar@<2ex>^{q^{2}_{v}} "11"; "12"};
{\ar@{.>}@<0ex>|{m^{2}_{v}} "11"; "12"};
{\ar@<-2ex>_{p^{2}_{v}} "11"; "12"};
%%%%%%%%%%%%%%%%%%%%%%%%%%%%%%%%%%%%%%%%%%%%%
{\ar@<2ex>^{c^{3}_{v}} "20"; "21"};
{\ar@<0ex>|{i^{3}_{v}} "21"; "20"};
{\ar@<-2ex>_{d^{3}_{v}} "20"; "21"};
{\ar@<2ex>^{q^{3}_{v}} "21"; "22"};
{\ar@{.>}@<0ex>|{m^{3}_{v}} "21"; "22"};
{\ar@<-2ex>_{p^{3}_{v}} "21"; "22"};
\endxy
\end{gathered}
\end{equation}
in which all rows and columns $A_{-,n}$ and $A_{m,-}$ are cocategories, and each trio of the form $f^{-}_{h}$ or $f^{-}_{v}$ is a morphism of cocategories.  We sometimes refer to the commutativity of the dotted square as the middle four interchange axiom in a double cocategory.\\
%Note that in the diagram defining a double cocategory the entire diagram is uniquely determined by its restriction to the three quadrants omitting $A_{3,3}$.  For $A_{3,3}=A_{2,3} +_{A_{1,3}}A_{2,3}=A_{3,2} +_{A_{3,1}}A_{3,2}$ and $f^{3}_{h}=f^{2}_{h}+_{{f^{1}_{h}}}f^{2}_{h}$ and $f^{3}_{v}=f^{2}_{v}+_{{f^{1}_{v}}}f^{2}_{v}$ for $f \in \{p,m,q\}$.  Indeed the only commutativity in the lower right quadrant not automatically implied the others is the square above right: this is the square expressing the \emph{middle four interchange} in a double cocategory, and is central to our investigations.  Given a diagram as above left in which each row and column is a cocategory and each commutativity required of a double cocategory holds \emph{excepting} the single square above right, we call the resulting structure a \emph{pre-double cocategory.}
Observe that if $\mathbf A$ and $\mathbf B$ are cocategories in $\f C$ and $\otimes:\f C \times \f C \to \f D$ is cocontinuous in each variable then the pointwise tensor product of $\mathbf A$ and $\mathbf B$ yields a double cocategory $\mathbf A \otimes \mathbf B$.  This follows from the fact that each of $A_{i} \otimes -$ and $- \otimes B_{j}$ preserves cocategories.  In particular given a bifunctor $\otimes:\Cat \times \Cat \to \f C$ cocontinuous in each variable the corresponding double cocategory in $\f C$ of Proposition~\ref{prop:equiv} is $\arr \otimes \arr$ as below:
\begin{equation}
\label{eq:tensor}
\begin{gathered}
\xy
(0,0)*+{1\otimes 1}="00"; (30,0)*+{\atwo \otimes 1}="10";(60,0)*+{\athree \otimes 1}="20";
%%%%%%%%%%%%%%%
(0,-30)*+{1 \otimes \atwo}="01"; (30,-30)*+{\atwo \otimes \atwo}="11";(60,-30)*+{\athree \otimes \atwo}="21";
%%%%%%%%%%%%%%%%
(0,-60)*+{1 \otimes \athree}="02"; (30,-60)*+{\atwo \otimes \athree}="12";(60,-60)*+{\athree \otimes \athree}="22";
%%%%%%%%%%%%%%%%%%%%%%%%%%%%%%%%%%%%%%%%%%%
{\ar@<2ex>^{d \otimes 1} "00"; "10"};
{\ar@<0ex>|{i \otimes 1} "10"; "00"};
{\ar@<-2ex>_{c \otimes 1} "00"; "10"};
{\ar@<2ex>^{p \otimes 1} "10"; "20"};
{\ar@<0ex>|{m \otimes 1} "10"; "20"};
{\ar@<-2ex>_{q \otimes 1} "10"; "20"};
%%%%%%%%%%%%%%%%%%%%%%%%%%%%%%%%%%%%%%%%%%%%%
{\ar@<2ex>^{d \otimes \atwo} "01"; "11"};
{\ar@<0ex>|{i \otimes \atwo} "11"; "01"};
{\ar@<-2ex>_{c \otimes \atwo} "01"; "11"};
{\ar@<2ex>^{p \otimes \atwo} "11"; "21"};
{\ar@<0ex>|{m \otimes \atwo} "11"; "21"};
{\ar@<-2ex>_{q \otimes \atwo} "11"; "21"};
%%%%%%%%%%%%%%%%%%%%%%%%%%%%%%%%%%%%%%%%%%%%%
{\ar@<2ex>^{d \otimes \athree} "02"; "12"};
{\ar@<0ex>|{i \otimes \athree} "12"; "02"};
{\ar@<-2ex>_{c \otimes \athree} "02"; "12"};
{\ar@<2ex>^{p \otimes \athree} "12"; "22"};
{\ar@<0ex>|{m \otimes \athree} "12"; "22"};
{\ar@<-2ex>_{q \otimes \athree} "12"; "22"};
%%%%%%%%%%%%%%%%%%%%%%%%%%%%%%%%%%%%%%%%%%%
{\ar@<2ex>^{1 \otimes c} "00"; "01"};
{\ar@<0ex>|{1 \otimes i} "00"; "01"};
{\ar@<-2ex>_{1 \otimes d} "00"; "01"};
{\ar@<2ex>^{1 \otimes q} "01"; "02"};
{\ar@<0ex>|{1 \otimes m} "01"; "02"};
{\ar@<-2ex>_{1 \otimes p} "01"; "02"};
%%%%%%%%%%%%%%%%%%%%%%%%%%%%%%%%%%%%%%%%%%%%
{\ar@<2ex>^{\atwo \otimes c} "10"; "11"};
{\ar@<0ex>|{\atwo \otimes i} "11"; "10"};
{\ar@<-2ex>_{\atwo \otimes d} "10"; "11"};
{\ar@<2ex>^{\atwo \otimes q} "11"; "12"};
{\ar@<0ex>|{\atwo \otimes m} "11"; "12"};
{\ar@<-2ex>_{\atwo \otimes p} "11"; "12"};
%%%%%%%%%%%%%%%%%%%%%%%%%%%%%%%%%%%%%%%%%%%%
{\ar@<2ex>^{\athree \otimes c} "20"; "21"};
{\ar@<0ex>|{\athree \otimes  i} "21"; "20"};
{\ar@<-2ex>_{\athree \otimes  d} "20"; "21"};
{\ar@<2ex>^{\athree \otimes q} "21"; "22"};
{\ar@<0ex>|{\athree \otimes m} "21"; "22"};
{\ar@<-2ex>_{\athree \otimes  p} "21"; "22"};
(68,-62)*+{.};
\endxy
\end{gathered}
\end{equation}
The cartesian product of categories $\atwo \times \atwo$ is of course the free commutative square
$$\xy
(0,0)*+{(0,0)}="00"; (15,0)*+{(1,0)}="10";(0,-15)*+{(0,1)}="01";(15,-15)*+{(1,1)}="11";
{\ar_{(0,f)} "00"; "01"};
{\ar^{(f,0)} "00"; "10"};
{\ar^{(1,f)} "10"; "11"};
{\ar_{(f,1)} "01"; "11"};
(7,-7)*+{=}
\endxy
\hspace{2cm}
\xy
(0,0)*+{(0,0)}="00"; (15,0)*+{(1,0)}="10";(0,-15)*+{(0,1)}="01";(15,-15)*+{(1,1)}="11";
{\ar_{(0,f)} "00"; "01"};
{\ar^{(f,0)} "00"; "10"};
{\ar^{(1,f)} "10"; "11"};
{\ar_{(f,1)} "01"; "11"};
(7,-7)*+{\neq}
\endxy$$ whilst the funny tensor product $\atwo \star \atwo$ is the non-commutative square.  These are depicted above.  The double cocategory structures $\arr \times \arr$ and $\arr \star \arr$ are clear in either case.
\subsection{Only two biclosed structures.}
In the case of both the funny and cartesian products, the associativity constraint $(A \otimes B) \otimes C \cong A \otimes (B \otimes C)$ is the only possible \emph{natural} isomorphism.  (This is clear at the triples $(1,1,1)$, $(\atwo,1,1), (1,\atwo,1)$ and $(1,1,\atwo)$, thus at all triples involving only $1$'s and $\atwo$'s by naturality, and these force the general case.)  Similar remarks apply to the unit isomorphisms.   To show that these are the only monoidal biclosed structures, it therefore suffices to prove:
\begin{Proposition}
The funny tensor product and cartesian product are the only two cocontinuous bifunctors with unit on \Cat.
\end{Proposition}
\noindent Any such bifunctor has unit $1$ by Proposition~\ref{prop:unit}.  Such bifunctors correspond to double cocategories as in \eqref{eq:tensor} in which the top and left cocategories $\arr \otimes 1$ and $1 \otimes \arr$ are isomorphic to $\arr$.  Accordingly, the above result will follow upon proving:
\begin{Theorem}\label{thm:CatDbl}
Up to isomorphism there exist just two double cocategories in \Cat whose top and left cocategories coincide as the arrow cocategory $\arr$ in \Cat.  They are $\arr \star \arr$ and $\arr \times \arr$.
\end{Theorem}
\noindent The following result from \cite{Lumsdaine2011A-small} is a helpful starting point.
\begin{Lemma}
Each cocategory in \Set is a co-equivalence relation and, consequently, the co-kernel pair of its equaliser.
\end{Lemma}
\noindent We sketch the proof.  To begin with, given a cocategory $\mathbf A$ in \Set as in \eqref{eq:cocat}, one should show that the maps $d,c:A_{1} \rightrightarrows A_{2}$ are jointly epi.   The pushout projections $p,q:A_{2} \rightrightarrows A_{3}$ certainly are, so given $x \in A_{2}$ there exists $y \in A_{2}$ such that $py=mx$ or $qy=mx$.  So either $x=m(i,1)x=dix$ or $x=m(1,i)x=cix$; therefore $d$ and $c$ are jointly epi and $\mathbf A$ a co-preorder.  We omit details of the symmetry.  Because $\Set^{op}$ is Barr-exact, being monadic over \Set by \cite{Pare1974Colimits}, each co-equivalence relation in \Set is the co-kernel pair of its equaliser.
\newline{}
Let us call $\Ob$ the cocategory in \Set:
$$\xy
(0,0)*+{1}="a0"; (40,0)*+{2}="b0";(80,0)*+{3}="c0";
{\ar@<1.5ex>^{d} "a0"; "b0"};
{\ar@<0ex>|{i} "b0"; "a0"};
{\ar@<-1.5ex>_{c} "a0"; "b0"};
{\ar@<1.5ex>^{p} "b0"; "c0"};
%{\ar@<1.5ex>|{r} "c0"; "b0"};
{\ar@<0ex>|{m} "b0"; "c0"};
%{\ar@<-1.5ex>|{l} "c0"; "b0"};
{\ar@<-1.5ex>_{q} "b0"; "c0"};
\endxy$$
which is the co-kernel pair of $\emptyset \to 1$, so that $d$ and $c$ are the coproduct inclusions.  Because colimits of cocategories are pointwise each set $X$ gives rise to a cocategory $X.\Ob$ by taking componentwise copowers (so $X.\Ob_{2}=X.2$ etc).  By exactness a morphism of cocategories $f:\mathbf A \to \Ob$ corresponds to a function $Eq(d_{A},c_{A}) \to \emptyset$ between equalisers, and a unique such exists just when $Eq(d_{A},c_{A})=\emptyset$.  Since $A_{1}$ is the co-kernel pair of its equaliser it follows that $A_{2}=A_{1}.2$ and that $f=!.\Ob$ for $!:A_{1} \to 1$.  Since each cocategory present in a double cocategory comes equipped with a map to either its top or left cocategory it follows that:
\begin{Corollary}\label{cor:SetDbl}
There exists only one double cocategory in \Set whose left and top cocategories coincide as the cocategory \Ob; namely the cartesian product $\Ob \times \Ob$.
\end{Corollary}
\begin{proofof}[Theorem~\ref{thm:CatDbl}]
By Proposition~\ref{prop:equiv} we are justified in denoting such a double cocategory $\arr \otimes \arr$ as in \eqref{eq:tensor}.
\begin{enumerate}
\item
Since the functor $(-)_{0}:\Cat \to \Set$ preserves colimits it preserves (double) cocategories, and takes the arrow cocategory $\arr$ to $\Ob$.   By Corollary~\ref{cor:SetDbl} we have an isomorphism of double cocategories $(\arr \otimes \arr)_{0} \cong \Ob \times \Ob$ and it is harmless to assume that they are equal.  In particular we then have $(\atwo \otimes \atwo)_{0}= 2 \times 2$.
Let us abbreviate $d \otimes \atwo$ by $d^{2}_{h}$ and so on, as in \eqref{eq:dbl}.  We then have a partial diagram of $\atwo \otimes \atwo$:
\begin{equation}
\label{eq:foursides}
\cd{(0,0) \ar[r]^{d^{2}_{v}f} \ar[d]_{d^{2}_{h}f} & (1,0)\ar[d]^{c^{2}_{h}f} \\
(0,1)\ar[r]_{c^{2}_{v}f} & (1,1)}
\end{equation}
which gives a complete picture on objects.  There may exist more arrows than depicted, but all morphisms are of the form $(a,b) \to (c,d)$ where $a \leq c$ and $b \leq d$.  If, for instance, there were an arrow $\alpha:d^{2}_{h}1=(0,1) \to (0,0)=d^{2}_{h}0$ we would obtain $i^{2}_{h}\alpha:1 \to 0 \in \atwo$ but no such arrow exists, and one rules out the other possibilities in a similar fashion.
\item
The next problem is to show that each of $d^{2}_{h}, c^{2}_{h}, d^{2}_{v}$ and $c^{2}_{v}$ is \emph{fully faithful}, which amounts to the assertion that there exists a unique arrow on each of the four sides of $\atwo \otimes \atwo$ (as in \eqref{eq:foursides}) and no non-identity endomorphisms.  The four cases are entirely similar and we will consider only $d^{2}_{h}$.  Certainly $d^{2}_{h}$ is faithful.  \emph{Suppose that $p^{2}_{h}$ were full when restricted to the full image of $d^{2}_{h}$}.  Then $\alpha:d^{2}_{h}x \to d^{2}_{h}y$ gives $m^{2}_{h}\alpha:m^{2}_{h}d^{2}_{h}x \to m^{2}_{h}d^{2}_{h}y$.  Since $m^{2}_{h}d^{2}_{h}=p^{2}_{h}d^{2}_{h}$ we obtain $\beta:d^{2}_{h}x \to d^{2}_{h}y$ such that $p^{2}_{h}\beta=m^{2}_{h}\alpha$.  But then $\alpha=(i^{2}_{h},1)m^{2}_{h}\alpha=(i^{2}_{h},1)p^{2}_{h}\beta=d^{2}_{h}i^{2}_{h}\beta$ so that $d^{2}_{h}$ is full as well as faithful.\newline{}
So it will suffice to show that $p^{2}_{h}$ is full on the objects $(0,0)$ and $(0,1)$ in the image of $d^{2}_{h}$.  This involves computing the pushout $\athree \otimes \atwo$ of categories and we can reduce this to a simpler computation involving graphs.  To this end observe that $\atwo=1 \otimes \atwo$ is free on the graph $\atwo_{G}=\{0 \to 1\}$.  By adjointness $d^{2}_{h}$ and $c^{2}_{h}$ correspond to maps $d^{\prime},c^{\prime}:\atwo_{G} \rightrightarrows U(\atwo \otimes \atwo)$ and we obtain a diagram:
$$\xy
(0,15)*+{F(\atwo_{G})}="a1"; (40,15)*+{FU(\atwo \otimes \atwo)}="b1";(80,15)*+{F(P)}="c1";
(0,0)*+{1 \otimes \atwo}="a0"; (40,0)*+{\atwo \otimes \atwo}="b0";(80,0)*+{\athree \otimes \atwo}="c0";
{\ar@<1ex>^{d^{2}_{h}} "a0"; "b0"};
{\ar@<-1ex>_{c^{2}_{h}} "a0"; "b0"};
{\ar@<1ex>^{p^{2}_{h}} "b0"; "c0"};
{\ar@<-1ex>_{q^{2}_{h}} "b0"; "c0"};
%%%%%%%%%%%%%%%%%%%%
{\ar@<1ex>^{Fd^{\prime}} "a1"; "b1"};
{\ar@<-1ex>_{Fc^{\prime}} "a1"; "b1"};
{\ar@<1ex>^{Fp^{\prime}} "b1"; "c1"};
{\ar@<-1ex>_{Fq^{\prime}} "b1"; "c1"};
%%%%%%%%%%%%%%%%%%%
{\ar@<0ex>_{1} "a1"; "a0"};
{\ar@<0ex>^{\epsilon} "b1"; "b0"};
{\ar@<0ex>^{k} "c1"; "c0"};
\endxy$$
where $\epsilon$ is an identity on objects and full functor, given by the counit of the adjunction $F \dashv U$.  $P$ is the pushout of $d^{\prime}$ and $c^{\prime}$ in the category of graphs so that the top row is a pushout of categories.  Since bijective-on-objects and full functors form the left class of a factorisation system on \Cat they are closed under pushout in $\Cat^{\atwo}$; thus the induced map $k$ between pushouts is bijective-on-objects and full too, and it is harmless to suppose that it is the identity on objects.  By a two from three argument it follows that $p^{2}_{h}$ will be full on $\{(0,0),(0,1)\}$ so long as $Fp^{\prime}$ is full on these same objects.  In the category of graphs the pushout $P$ of the monos $d^{\prime}$ and $c^{\prime}$ is easily calculated, and it is clear that $p^{\prime}$ is a full embedding of graphs where restricted to $(0,0)$ and $(0,1)$. $P$ has the further property that the only morphisms in $P$ having codomain $(0,0)$ or $(0,1)$ also have domain amongst these two objects.  Using the explicit construction of morphisms in $FP$ as paths in $P$, these two facts are enough to ensure that $Fp^{\prime}$ is full where restricted to $(0,0)$ and $(0,1)$, as required.
\item
Because each of $d^{2}_{h}, c^{2}_{h}, d^{2}_{v}$ and $c^{2}_{v}$ is fully faithful our knowledge of $\atwo \otimes \atwo$ is complete with the exception that we have not determined the hom-set $\atwo \otimes \atwo((0,0),(1,1))$.    The two cases $\atwo \otimes \atwo((0,0),(1,1))=1,2$ correspond to the free commutative square and non-commutative square and these do give double cocategories.\newline{}
To rule out any other cases suppose that the complement $$X=\atwo \otimes \atwo((0,0),(1,1)) \setminus \{c^{2}_{h}f \circ d^{2}_{v}f, c^{2}_{v}f \circ d^{2}_{h}f\}$$ is non-empty.  Then the pushout $\atwo \otimes \athree$ is generated by the graph below right
$$\cd{(0,0) \ar[r]^{d^{2}_{v}f} \ar[dr]^{X} \ar[d]_{d^{2}_{h}f} & (1,0)\ar[d]^{c^{2}_{h}f} & &  (0,0) \ar[r]^{p^{2}_{h}d^{2}_{v}f} \ar[dr]^{X_{l}} \ar[d]_{p^{2}_{h}d^{2}_{h}f} & (1,0)\ar[dr]^{X_{r}} \ar[r]^{q^{2}_{h}d^{2}_{v}f}\ar[d]|{p^{2}_{h}c^{2}_{h}f} & (2,0)\ar[d]^{q^{2}_{h}c^{2}_{h}f}\\
(0,1)\ar[r]_{c^{2}_{v}f} & (1,1) & & (0,1)\ar[r]_{p^{2}_{h}c^{2}_{v}f} & (1,1) \ar[r]_{q^{2}_{h}c^{2}_{h}f} & (2,1) }$$
 subject to the equations asserting the commutativity of the left and right squares whenever the square does in fact commute in $\atwo \otimes \atwo$.  We have two copies $X_{l}$ and $X_{r}$ of $X$ as depicted.  At $\alpha \in X$ we must have $m^{2}_{h}\alpha:(0,0) \to (2,1)$ and also $(i^{2}_{h},1)m^{2}_{h}\alpha=(1,i^{2}_{h})m^{2}_{h}\alpha=\alpha$.  This implies that the path $m^{2}_{h}\alpha$ must involve both a morphism of $X_{l}$ and a morphism of $X_{r}$, but there exists no such path.
\end{enumerate}
\end{proofof}

\section{Monoidal biclosed structures on \twocat}
There are at least five monoidal biclosed structures on \twocat.  %Of these, three are symmetric monoidal and  only one is a monoidal model structure.
By Remark~\ref{rk:HigherUnit} any monoidal biclosed structure on \twocat has unit $1$ and, correspondingly, the objects of any possible hom 2-category must be the 2-functors.  The 1-cells between these are, in the five cases, 2-natural, lax natural, pseudonatural, oplax natural and not-necessarily-natural transformations; in each case, modifications are the 2-cells.
$$
\cd{& Lax(A,B)\ar[dr] & & & & A \otimes_{l} B \ar[dr] \\
[A,B]\ar[ur] \ar[r] \ar[dr] & Ps(A,B) \ar[r] &[A,B]_{f} & & A \star B \ar[ur] \ar[r] \ar[dr] & A \otimes_{p} B \ar[r] & A \times B\\
& Oplax(A,B)\ar[ur] & & & & A \otimes_{o} B \ar[ur] \rlap{   .}}
$$
The five homs are connected by evident inclusions and forgetful functors as in the diagram above left.  The arrows therein are the comparison maps of closed functor structures on the identity $1:\twocat \to \twocat$.  By adjointness (of the $A \otimes - \dashv [A,-]$ variety) there is a corresponding diagram of functors on the right, each of which is now the component of an opmonoidal structure on $1:\twocat \to \twocat$.  Reading left to right and top to bottom we have the funny tensor products, lax/ordinary and oplax Gray tensor products \cite{Gray1974Formal}, and the cartesian product.\newline{}
The tensor product $\otimes_{l}$ is monoidal biclosed; indeed $$\twocat(A,Lax(B,C)) \cong \twocat(A \otimes_{l} B,C) \cong \twocat(B,Oplax(A,C)) \hspace{0.2cm}\textnormal{.}$$  Likewise $\otimes_{o}$ is biclosed, but neither tensor product is symmetric; rather, they are opposite: $A \otimes_{l}B \cong B \otimes_{o} A$.  The other cases are symmetric monoidal.\newline{}
Viewing $\atwo$ as a 2-category in which each 2-cell is an identity, let us calculate $\atwo \otimes_{l} \atwo$ and $\atwo \otimes_{p} \atwo$.  Both have underlying category the non-commuting square $\atwo \star \atwo$ and are depicted, in turn, below.
\begin{equation}\label{eq:square}
\begin{gathered}
\xy
(0,0)*+{(0,0)}="00"; (15,0)*+{(1,0)}="10";
(0,-15)*+{(0,1)}="01";(15,-15)*+{(1,1)}="11";
{\ar^{(f,0)} "00"; "10"};
{\ar_{(0,f)} "00"; "01"};
{\ar^{(1,f)} "10"; "11"};
{\ar_{(f,1)} "01"; "11"};
{\ar@{=>}^{}(11,-5)*+{};(5,-10)*+{}};
\endxy
\hspace{1cm}
\xy
(0,0)*+{(0,0)}="00"; (15,0)*+{(1,0)}="10";
(0,-15)*+{(0,1)}="01";(15,-15)*+{(1,1)}="11";
{\ar^{(f,0)} "00"; "10"};
{\ar_{(0,f)} "00"; "01"};
{\ar^{(1,f)} "10"; "11"};
{\ar_{(f,1)} "01"; "11"};
(8,-7)*+{\cong};
\endxy
\end{gathered}
\end{equation}
We do not know whether the above five encompass all monoidal biclosed structures on \twocat, but can say that the techniques of Section 2 are insufficient to determine whether this is true.  Those techniques worked because each cocontinuous bifunctor with unit on \Cat underlies a \emph{unique} monoidal structure, whereas on \twocat there exist such bifunctors not underlying any monoidal structure.  For an example let $CB_{1}(\twocat)$ denote the category of cocontinuous bifunctors on \twocat with unit $1$.  It is easily seen that the forgetful functor $CB_{1}(\twocat) \to [\twocat^{2},\twocat]$ creates connected colimits.  In particular the co-kernel pair of $\star \to \otimes_{l}$ yields an object $\otimes_{2}$ of $CB_{1}(\twocat)$, with $\atwo \otimes_{2} \atwo$ given by
$$\xy
(0,0)*+{(0,0)}="00"; (15,0)*+{(1,0)}="10";
(0,-15)*+{(0,1)}="01";(15,-15)*+{(1,1)}="11";
{\ar^{(f,0)} "00"; "10"};
{\ar_{(0,f)} "00"; "01"};
{\ar^{(1,f)} "10"; "11"};
{\ar_{(f,1)} "01"; "11"};
{\ar@{=>}^{}(11,-3)*+{};(5,-7)*+{}};
{\ar@{=>}^{}(11,-7)*+{};(5,-11)*+{}};
(25,-16)*+{.};
\endxy$$
One can show that this tensor product is not associative up to natural isomorphism as follows.  Any natural isomorphism would have component at $(\atwo \otimes _{2} \atwo) \otimes_{2} \atwo \to \atwo \otimes_{2} \atwo (\otimes_{2} \atwo)$ the associativity isomorphism $(\atwo \star \atwo) \star \atwo \to \atwo \star (\atwo \star \atwo)$ on underlying categories, but direct calculation shows that this last isomorphism cannot be extended to 2-cells.
\section{The case of Gray-categories}
%In \twocat there are infinitely many double cocategories with top and left cocategory $\arr$.
The inclusion $j:\Cat \to \twocat$ that views each category as a locally discrete 2-category allows us to view the cocategory $\arr$ of \eqref{eq:arr} as a cocategory in $\twocat$, so that we obtain the tensor product double cocategories $\arr \otimes_{l} \arr$ and $\arr\otimes_{p} \arr$ corresponding to the lax and ordinary Gray tensor product.  %The forgetful functor $(-)_{1}:\twocat \to \Cat$ is strong monoidal from both $\otimes_{l}$ and $\otimes_{p}$ to $\star$ so that in each case $(\arr \otimes \arr)_{1}=\arr \star \arr$.
The 2-categories $\atwo \otimes_{l} \atwo$ and $\atwo \otimes_{p} \atwo$, depicted in \eqref{eq:square}, each consist of a non-commuting square and 2-cell -- invertible in the second case -- within.  Homming from $\arr \otimes_{l} \arr$ into a 2-category \f C yields the double category $\twocat(S \otimes_{l} S,\f C)$ of quintets in \f C \cite{Ehresmann1963Categories}, whose \emph{squares} are lax squares in $\f C$ as below left.  The fact that the middle four interchange axiom holds in the double cocategory $\arr \otimes_{l}\arr$ corresponds, by Yoneda, to the fact that the two ways of composing a diagram of 2-cells as below right

\begin{equation}
\label{eq:config}
\begin{gathered}
\xy
(0,-7.5)*+{.}="00"; (15,-7.5)*+{.}="10";
(0,-22.5)*+{.}="01";(15,-22.5)*+{.}="11";
{\ar^{} "00"; "10"};
{\ar_{} "00"; "01"};
{\ar^{} "10"; "11"};
{\ar_{} "01"; "11"};
{\ar@{=>}^{}(11,-12.5)*+{};(5,-17.5)*+{}};
\endxy
\hspace{2cm}
\xy
(0,0)*+{.}="00"; (15,0)*+{.}="10"; (30,0)*+{.}="20";
(0,-15)*+{.}="01";(15,-15)*+{.}="11"; (30,-15)*+{.}="21";
(0,-30)*+{.}="02";(15,-30)*+{.}="12"; (30,-30)*+{.}="22";
{\ar^{} "00"; "10"}; {\ar^{} "10"; "20"};
{\ar_{} "01"; "11"}; {\ar_{} "11"; "21"};
{\ar_{} "02"; "12"}; {\ar_{} "12"; "22"};
{\ar_{} "00"; "01"}; {\ar^{} "01"; "02"};
{\ar_{} "10"; "11"}; {\ar^{} "11"; "12"};
{\ar_{} "20"; "21"}; {\ar^{} "21"; "22"};
{\ar@{=>}^{}(11,-5)*+{};(5,-10)*+{}};{\ar@{=>}^{}(26,-5)*+{};(20,-10)*+{}};
{\ar@{=>}^{}(11,-20)*+{};(5,-25)*+{}};{\ar@{=>}^{}(26,-20)*+{};(20,-25)*+{}};
\endxy
\end{gathered}
\end{equation}
in a 2-category always coincide.\newline{}
In a Gray-category it is not necessarily the case that the two composite 2-cells coincide and, correspondingly, neither $\arr \otimes_{l} \arr$ nor $\arr \otimes_{p} \arr$ remains a double cocategory upon taking its image under the inclusion $j:\twocat \to \GCat$.\begin{footnote}{Specifically, the inclusion fails to preserve the pushouts $\athree \otimes_{l} \athree$ and $\athree \otimes_{p} \athree$.}\end{footnote}  This lack of double cocategories in $\GCat$ ultimately limits the biclosed structures that can exist on $\GCat$ -- see Theorem~\ref{thm:TensorGray} below.
\subsection{Gray-categories and sesquicategories.}
Since the above problem concerns the failure of the middle four interchange for $2$-cells, it manifests itself at the simpler two-dimensional level of \emph{sesquicategories} \cite{Street1996Categorical}.    Let us recall Gray-categories, sesquicategories and their inter-relationship.\newline{}
A \emph{Gray-category} is a category enriched in $(\twocat,\otimes_{p})$ and we write $\GCat$ for the category of Gray-categories.  A sesquicategory $\f C$ is a category enriched in $(\Cat,\star)$.  In elementary terms such a \f C consists of a set of objects; for each pair $B,C$ of objects a hom-category $\f C(B,C)$ of 1-cells and 2-cells, together with whiskering functors $\f C(f,C):\f C(B,C) \to \f C(A,C)$ and $\f C(B,g):\f C(B,C) \to \f C(B,D)$ for $f:A \to B$ and $g:C \to D$ such that $(g\alpha)f=g(\alpha f)$ where defined. We write \SCat for the category of sesquicategories.\newline{}
The forgetful functor $(-)_{1}:\twocat \to \Cat$ has both left and right adjoints.  Furthermore $(-)_{1}:(\twocat,\otimes_{p}) \to (\Cat,\star)$ is strong monoidal and is consequently the left adjoint of a monoidal adjunction \cite{Kelly1974Doctrinal}.  It follows that the lifted functor $(-)_{2}:\GCat \to \SCat$, which forgets 3-cells, is itself a left adjoint.  This was observed in \cite{Lack2011A-quillen}.  In particular $(-)_{2}$ preserves colimits and so double cocategories.  We have a commutative triangle of truncation functors: $$\xy
(0,0)*+{\GCat}="00"; (40,0)*+{\SCat}="10";
(20,-15)*+{\Cat}="11";
{\ar^{(-)_{2}} "00"; "10"};
{\ar^{(-)_{1}} "00"; "11"};
{\ar_{(-)_{1}} "10"; "11"};
{\ar@<2ex>^{j} "11"; "00"};
{\ar@<-2ex>_{j} "11"; "10"};
%(45,-15)*+{\rlap{.}};
\endxy$$
in which both functors to $\Cat$ have fully faithful left adjoints.  Denoted by $j$ in either case, these adjoints view a category as a Gray-category or sesquicategory whose higher dimensional cells are identities.\newline{}
In what follows the term \emph{pre-double cocategory} will refer to a diagram as below left:
 $$\xy
(0,0)*+{A_{1,1}}="00"; (25,0)*+{A_{2,1}}="10";(50,0)*+{A_{3,1}}="20";
%%%%%%%%%%%%%%%
(0,-25)*+{A_{1,2}}="01"; (25,-25)*+{A_{2,2}}="11";(50,-25)*+{A_{3,2}}="21";
%%%%%%%%%%%%%%%%
(0,-50)*+{A_{1,3}}="02"; (25,-50)*+{A_{2,3}}="12";%(50,-50)*+{A_{3,3}}="22";
%%%%%%%%%%%%%%%%%%%%%%%%%%%%%%%%%%%%%%%%%%%
{\ar@<2ex>^{d^{1}_{h}} "00"; "10"};
{\ar@<0ex>|{i^{1}_{h}} "10"; "00"};
{\ar@<-2ex>_{c^{1}_{h}} "00"; "10"};
{\ar@<2ex>^{p^{1}_{h}} "10"; "20"};
{\ar@<0ex>|{m^{1}_{h}} "10"; "20"};
{\ar@<-2ex>_{q^{1}_{h}} "10"; "20"};
%%%%%%%%%%%%%%%%%%%%%%%%%%%%%%%%%%%%%%%%%%%%%
{\ar@<2ex>^{d^{2}_{h}} "01"; "11"};
{\ar@<0ex>|{i^{2}_{h}} "11"; "01"};
{\ar@<-2ex>_{c^{2}_{h}} "01"; "11"};
{\ar@<2ex>^{p ^{2}_{h}} "11"; "21"};
{\ar@<0ex>|{m^{2}_{h}} "11"; "21"};
{\ar@<-2ex>_{q^{2}_{h}} "11"; "21"};
%%%%%%%%%%%%%%%%%%%%%%%%%%%%%%%%%%%%%%%%%%%%%
{\ar@<2ex>^{d^{3}_{h}} "02"; "12"};
{\ar@<0ex>|{i^{3}_{h}} "12"; "02"};
{\ar@<-2ex>_{c^{3}_{h}} "02"; "12"};
%{\ar@{.>}@<2ex>^{p ^{3}_{h}} "12"; "22"};
%{\ar@{.>}@<0ex>|{m^{3}_{h}} "12"; "22"};
%{\ar@{.>}@<-2ex>_{q^{3}_{h}} "12"; "22"};
%%%%%%%%%%%%%%%%%%%%%%%%%%%%%%%%%%%%%%%%%%%%%
{\ar@<2ex>^{c^{1}_{v}} "00"; "01"};
{\ar@<0ex>|{i^{1}_{v}} "01"; "00"};
{\ar@<-2ex>_{d^{1}_{v}} "00"; "01"};
{\ar@<2ex>^{q^{1}_{v}} "01"; "02"};
{\ar@<0ex>|{m^{1}_{v}} "01"; "02"};
{\ar@<-2ex>_{p^{1}_{v}} "01"; "02"};
%%%%%%%%%%%%%%%%%%%%%%%%%%%%%%%%%%%%%%%%%%%
{\ar@<2ex>^{c^{2}_{v}} "10"; "11"};
{\ar@<0ex>|{i^{2}_{v}} "11"; "10"};
{\ar@<-2ex>_{d^{2}_{v}} "10"; "11"};
{\ar@<2ex>^{q^{2}_{v}} "11"; "12"};
{\ar@<0ex>|{m^{2}_{v}} "11"; "12"};
{\ar@<-2ex>_{p^{2}_{v}} "11"; "12"};
%%%%%%%%%%%%%%%%%%%%%%%%%%%%%%%%%%%%%%%%%%%%%
{\ar@<2ex>^{c^{3}_{v}} "20"; "21"};
{\ar@<0ex>|{i^{3}_{v}} "21"; "20"};
{\ar@<-2ex>_{d^{3}_{v}} "20"; "21"};
%{\ar@{.>}@<2ex>^{q^{3}_{v}} "21"; "22"};
%{\ar@{.>}@<0ex>|{m^{3}_{v}} "21"; "22"};
%{\ar@{.>}@<-2ex>_{p^{3}_{v}} "21"; "22"};
\endxy
\hspace{1cm}
\xy
(0,-12)*+{A_{2,2}}="11"; (25,-12)*+{A_{3,2}}="21";
(0,-37)*+{A_{2,3}}="12"; (25,-37)*+{A_{3,3}}="22";
%%%%%%%%%%%%%%%%%%%%%%%%%
{\ar@<2ex>^{q^{3}_{v}} "21"; "22"};
{\ar@{.>}@<0ex>|{m^{3}_{v}} "21"; "22"};
{\ar@<-2ex>_{p^{3}_{v}} "21"; "22"};
%%%%%%%%%%%%%%%%%%%%%%%%%
{\ar@<2ex>^{q^{2}_{v}} "11"; "12"};
{\ar@{.>}@<0ex>|{m^{2}_{v}} "11"; "12"};
{\ar@<-2ex>_{p^{2}_{v}} "11"; "12"};
%%%%%%%%%%%%%%%%%%%%%%%%%
{\ar@<2ex>^{p ^{3}_{h}} "12"; "22"};
{\ar@{.>}@<0ex>|{m^{3}_{h}} "12"; "22"};
{\ar@<-2ex>_{q^{3}_{h}} "12"; "22"};
%%%%%%%%%%%%%%%%%%%%%%%%%
{\ar@<2ex>^{p ^{2}_{h}} "11"; "21"};
{\ar@{.>}@<0ex>|{m^{2}_{h}} "11"; "21"};
{\ar@<-2ex>_{q^{2}_{h}} "11"; "21"};
\endxy$$
in which all columns $A_{m,-}$ and rows $A_{-,n}$ of length three are cocategories, and each trio of the form $f^{-}_{h}$ or $f^{-}_{v}$ is a morphism of cocategories.\newline{}
Being a double cocategory is merely a property of a pre-double cocategory: $A_{3,3}$ is uniquely determined as the pushout $A_{2,3} +_{A_{1,3}}A_{2,3}=A_{3,2} +_{A_{3,1}}A_{3,2}$ and $f^{3}_{h}=f^{2}_{h}+_{{f^{1}_{h}}}f^{2}_{h}$ and $f^{3}_{v}=f^{2}_{v}+_{{f^{1}_{v}}}f^{2}_{v}$ for $f \in \{p,m,q\}$. In fact, the only distinction is that in a pre-double cocategory the dotted square expressing \emph{middle four interchange} need not commute.\newline{}
The important example for us, $\mathbf I$, is a pre-double cocategory in $\SCat$.  Its image under $(-)_{1}:\SCat \to \Cat$ is $\arr \star \arr$ whilst each of its component sesquicategories $I_{n,m}$ is locally indiscrete, i.e., there is a unique invertible 2-cell between each parallel pair of 1-cells.  This suffices for a full description, but for the reader's convenience we point out that $I_{2,2}$ is the pseudo-commutative square $\atwo \otimes_{p} \atwo$ and indeed that $\mathbf I$ is the image of the pre-double cocategory underlying $\arr \otimes_{p} \arr$ under the inclusion $\twocat \to \SCat$.\newline{}
In order to see that $\mathbf I$ is not a double cocategory in \SCat it suffices, by Yoneda, to exhibit a configuration of shape \eqref{eq:config} of invertible 2-cells in a sesquicategory for which the two possible composites do not coincide.  This is straightforward.

\begin{Proposition}\label{prop:sesquicat}
There exist only two double cocategories in \SCat whose top and left cocategories coincide as the arrow cocategory $\arr$.  They are the double cocategories $\arr \star \arr$ and $\arr \times \arr$ in \Cat viewed as double cocategories in \SCat.
\end{Proposition}
\pf
By Proposition~\ref{prop:equiv} each double cocategory in \SCat is the tensor double cocategory $\arr \otimes \arr$ associated to a cocontinuous bifunctor $\otimes:\Cat \times \Cat \to \SCat$, and accordingly we will use the tensor notation of \eqref{eq:tensor}. %QQ Should this double cocat not be S \otimes T in general, but then our hypothesis is that S=T?
\newline{}
The forgetful functor $(-)_{1}:\SCat \to \Cat$ is cocontinuous and so preserves double cocategories.  It also preserves $\arr$ and so, by Theorem~\ref{thm:CatDbl}, either $(\arr \otimes \arr)_{1} = \arr \star \arr$ or $(\arr \otimes \arr)_{1} = \arr \times  \arr$.  In particular  $(\atwo \otimes \atwo)_{1}=\atwo \star \atwo$ or $(\atwo \otimes \atwo)_{1}=\atwo \times \atwo$.  Our main task is to show that all 2-cells in $\atwo \otimes \atwo$ are identities and we begin by treating the two cases together.\newline{}
\emph{Let us firstly show that each of $d^{2}_{h}, c^{2}_{h}, d^{2}_{v}$ and $c^{2}_{v}$ is a full embedding of sesquicategories}: thus all 2-cells on the four sides of $\atwo \otimes \atwo$ are identities.  The argument here is a straightforward adaptation of the second part of the proof of Theorem~\ref{thm:CatDbl} and we only outline it.  As there, it suffices to show that $d^{2}_{h}$ is fully faithful (now on 2-cells as well as 1-cells) and for this it is enough to show that $p^{2}_{h}$ is full when restricted to the full image of $d^{2}_{h}$.    In the \Cat-case of Theorem~\ref{thm:CatDbl} this amounted to studying a pushout of categories.  This was simplified to computing a pushout of graphs by using the adjunction $F \dashv U$ between categories and graphs, that the counit of this adjunction is bijective on objects and full, and that such functors are closed under pushouts in $\Cat^{\atwo}$.  In the sesquicategory case, we argue in a similar fashion, but replace the category of graphs by the category $\Der$ of \emph{derivation schemes} \cite{Street1996Categorical}: such are simply categories $\f A$ equipped with, for each parallel pair $f,g:A \to B \in \f A$, a set $\f A(A,B)(f,g)$ of 2-cells.  There is an evident forgetful functor $V:\SCat \to \Der$ which has a left adjoint $G=G_{2}G_{1}$.  The first component $G_{1}$ adds formally whiskered 2-cells $f\alpha g$ for triples appropriately aligned -- as in $f:A \to B$, $\alpha \in \f A(B,C)(h,k)$ and $g:C \to D$ -- whilst $G_{2}$ applies the free category construction to each hom graph $G_{1}\f A(A,B)$.  Each component of the counit $GV \to 1$ is an isomorphism on underlying categories and locally full.  Such sesquifunctors form the left class of a factorisation system on \SCat and are therefore closed under pushout.  By employing these facts as in Theorem~\ref{thm:CatDbl} the problem reduces to calculating a simple pushout in $\Der$.
\newline{}
\emph{Now suppose $(\atwo \otimes \atwo)_{1}=\atwo \star \atwo$}.  Consider the pre-double cocategory $\mathbf I$ discussed above in which $I_{2,2}=\atwo \otimes_{p} \atwo$ is the free pseudo-commutative square and in which both of $I_{3,2}$ and $I_{2,3}$ are locally indiscrete.  Indiscreteness induces a unique morphism $F:\arr \otimes \arr \to \mathbf I$ of pre-double cocategories whose image under $(-)_{1}$ is the identity.  The fact that $F$ is a morphism of pre-double cocategories implies that of the six faces of the cube:

\begin{equation}
\begin{gathered}
\xy
(10,10)*+{\atwo \otimes \athree}="11"; (30,10)*+{I_{2,3}}="31";
(0,00)*+{\atwo \otimes \atwo}="00";(20,0)*+{\athree \otimes \athree}="20";(40,0)*+{I_{3,3}}="40";
(10,-10)*+{\athree \otimes \atwo}="1-1";(30,-10)*+{I_{3,2}}="3-1";
{\ar^{F_{2,3}} "11"; "31"};
{\ar_{F_{3,2}} "1-1"; "3-1"};
{\ar^{\atwo \otimes m} "00"; "11"};
{\ar_{m \otimes \atwo} "00"; "1-1"};
{\ar^{m \otimes \athree} "11"; "20"};
{\ar_{\athree \otimes m} "1-1"; "20"};
{\ar^{m^{3}_{h}} "31"; "40"};
{\ar_{m^{3}_{v}} "3-1"; "40"};
{\ar^{F_{3,3}} "20"; "40"};
\endxy
\hspace{1.5cm}
\xy
(10,10)*+{\atwo \otimes \athree}="11"; (30,10)*+{I_{2,3}}="31";
(0,00)*+{\atwo \otimes \atwo}="00";(20,0)*+{I_{2,2}}="20";(40,0)*+{I_{3,3}}="40";
(10,-10)*+{\athree \otimes \atwo}="1-1";(30,-10)*+{I_{3,2}}="3-1";
{\ar^{F_{2,3}} "11"; "31"};
{\ar_{F_{3,2}} "1-1"; "3-1"};
{\ar^{\atwo \otimes m} "00"; "11"};
{\ar_{m \otimes \atwo} "00"; "1-1"};
{\ar^{F_{2,2}} "00"; "20"};
{\ar^{m^{3}_{h}} "31"; "40"};
{\ar_{m^{3}_{v}} "3-1"; "40"};
{\ar^{m^{2}_{v}} "20"; "31"};
{\ar_{m^{2}_{h}} "20"; "3-1"};
\endxy
\end{gathered}
\end{equation}
the five left-most are commutative.  In particular, the outer paths of the right diagram commute as do its two leftmost squares.  If in $\atwo \otimes \atwo$ there existed a 2-cell as below (or in the opposite direction)
$$\xy
(0,0)*+{(0,0)}="00"; (15,0)*+{(1,0)}="10";
(0,-15)*+{(0,1)}="01";(15,-15)*+{(1,1)}="11";
{\ar^{(f,0)} "00"; "10"};
{\ar_{(0,f)} "00"; "01"};
{\ar^{(1,f)} "10"; "11"};
{\ar_{(f,1)} "01"; "11"};
{\ar@{=>}^{}(11,-5)*+{};(5,-10)*+{}};
\endxy$$
then $F_{2,2}:\atwo \otimes \atwo \to I_{2,2}$ would be epi, and it would follow that the rightmost square of (4.2) commutes. Since $\mathbf I$ is not a double cocategory it does not commute and accordingly there exists no such 2-cell.\newline{}
Consequently there may only exist endo 2-cells of the two paths $(0,0) \rightrightarrows (1,1)$ of $\atwo \otimes \atwo$.  Now at the level of underlying categories the pushout is freely generated by the graph below right  $$\cd{(0,0) \ar[r]^{(f,0)} \ar[d]_{(0,f)} & (1,0)\ar[d]^{(1,f)} & &  (0,0) \ar[r]^{(g,0)}  \ar[d]_{(0,f)} & (1,0) \ar[r]^{(h,0)}\ar[d]|{(1,f)} & (2,0)\ar[d]^{(2,f)}\\
(0,1)\ar[r]_{(f,1)} & (1,1) & & (0,1)\ar[r]_{(g,1)} & (1,1) \ar[r]_{(h,1)} & (2,1) }$$
with the maps $p^{2}_{h}$ and $q^{2}_{h}$ including the left and right squares respectively.  It follows that the only 2-cells in $\athree \otimes \atwo$ are of the form $(h,1)(p^{2}_{h}\theta)$ and $(q^{2}_{h}\phi)(g,0)$.  Suppose that an endo 2-cell $\theta$ of $(1,f)(f,0)$ exists in $\atwo \otimes \atwo$.  It follows that $m^{2}_{h}\theta$ is an endomorphism of $(2,f)(h,0)(g,0)$ and so of the form $(q^{2}_{h}\phi)(g,0)$.  But then $\theta = (1,i^{2}_{h})m^{2}_{h}\theta = (1,i^{2}_{h})(q^{2}_{h}\theta) \circ (1,i^{2}_{h})(g,0)$ which is an identity 2-cell since  $(1,i^{2}_{h})q^{2}_{h}\theta = c^{2}_{h}i^{2}_{h}\theta$ is one.  Similarly each endomorphism of $(f,1)(0,f)$ is an identity and we finally conclude that all 2-cells in $\atwo \otimes \atwo$ are identities.\newline{}
\emph{Now suppose that $(\atwo \otimes \atwo)_{1}=\atwo \times \atwo$}.   As there exists only a single arrow $(0,0) \rightarrow (1,1)$ the hom-category $A=\atwo \otimes \atwo((0,0),(1,1))$ is a monoid, and we must prove that it is the trivial monoid.  By comparing the universal property of the pushout $\athree \otimes \atwo$ with that of the coproduct of monoids, we easily obtain the following complete description of $\athree \otimes \atwo$.  It has underlying category
% It remains to show that the monoid $A=\atwo \otimes \atwo((0,0),(1,1))$ is trivial; note that there is only a single arrow $(0,0) \rightarrow (1,1)$ by assumption, so $A$ is indeed a monoid. The pushout $\athree \otimes \atwo$ has underlying category
the product $\athree \times \atwo$ and the pushout inclusions $p^{2}_{h},q^{2}_{h}:\atwo \otimes \atwo \rightrightarrows \athree \otimes \atwo$ are full embeddings; thus $\athree \otimes \atwo((0,0),(1,1))=A=\athree \otimes \atwo((1,0),(2,1))$ as depicted left below.
$$\xy
(0,-7)*+{(0,0)}="00"; (15,-7)*+{(1,0)}="10";(30,-7)*+{(2,0)}="20";
(0,-22)*+{(0,1)}="01";(15,-22)*+{(1,1)}="11";(30,-22)*+{(2,1)}="21";
{\ar^{(g,0)} "00"; "10"};
{\ar_{(0,f)} "00"; "01"};
{\ar|{(1,f)} "10"; "11"};
{\ar^{(h,0)} "10"; "20"};
{\ar^{(2,f)} "20"; "21"};
{\ar_{(g,1)} "01"; "11"};
{\ar_{(h,1)} "11"; "21"};
(7,-15)*+{A};
(22,-15)*+{A};
\endxy
\hspace{0.5cm}
\xy
(0,0)*+{(0,0)}="00"; (15,0)*+{(1,0)}="10";
(0,-15)*+{(0,1)}="01";(15,-15)*+{(1,1)}="11";
(0,-30)*+{(0,2)}="02";(15,-30)*+{(1,2)}="12";
{\ar^{(f,0)} "00"; "10"};
{\ar_{(f,1)} "01"; "11"};
{\ar_{(0,g)} "00"; "01"};
{\ar^{(1,g)} "10"; "11"};
{\ar_{(0,h)} "01"; "02"};
{\ar^{(1,h)} "11"; "12"};
{\ar_{(f,2)} "02"; "12"};
(7,-23)*+{A};
(7,-8)*+{A};
\endxy
\hspace{0.5cm}
\xy
(0,0)*+{(0,0)}="00"; (15,0)*+{(1,0)}="10";(30,0)*+{(2,0)}="20";
(0,-15)*+{(0,1)}="01";(15,-15)*+{(1,1)}="11";(30,-15)*+{(2,1)}="21";
(0,-30)*+{(0,2)}="02";(15,-30)*+{(1,2)}="12";(30,-30)*+{(2,2)}="22";
{\ar^{(g,0)} "00"; "10"};
{\ar_{(0,g)} "00"; "01"};
{\ar|{(1,g)} "10"; "11"};
{\ar^{(h,0)} "10"; "20"};
{\ar^{(2,g)} "20"; "21"};
{\ar_{(0,h)} "01"; "02"};
{\ar|{(1,h)} "11"; "12"};
{\ar^{(2,h)} "21"; "22"};
{\ar_{(g,2)} "02"; "12"};
{\ar_{(h,2)} "12"; "22"};
{\ar_{(g,1)} "01"; "11"};
{\ar_{(h,1)} "11"; "21"};
(7,-8)*+{A};
(22,-8)*+{A};
(7,-23)*+{A};
(22,-23)*+{A};
\endxy$$
Moreover $\athree \otimes \atwo((0,0),(2,1))$ is the coproduct of monoids $$\athree \otimes \atwo((0,0),(1,1)) + \athree \otimes \atwo((1,0),(2,1)) = A+A$$ with inclusions given by whiskering.  The isomorphic $\atwo \otimes \athree$ is depicted centre above, also with hom monoid $\atwo \otimes \athree((0,0),(1,2))$ the coproduct $A+A$.  Likewise we calculate that $\athree \otimes \atwo$ and $\atwo \otimes \athree$ embed fully into the pushout $\athree \otimes \athree$ which is depicted above right, and that the monoid $\athree \otimes \athree((0,0)(2,2))$ is the coproduct $4.A$ with the coproduct inclusions given by whiskering as before.  %(All of these claims can be verified simply by comparing the universal properties of free monoids and of the defining pushouts of sesquicategories.)
\newline{}
As, for instance, discussed in Section 9.6 of \cite{Bergman1998An-invitation}, elements of a coproduct of monoids $\Sigma_{i \in I} B_{i}$ admit a \emph{unique normal form}: as words $b_{1}\ldots b_{n}$ where each $b_{i}$ belongs to some $B_{j} - \{e\}$ and no adjacent pair $(b_{i},b_{i+1})$ belong to the same $B_{j}$.  In order for this formulation to make sense the monoids $B_{j}$ must have disjoint underlying sets, or else must be replaced by disjoint isomorphs.  In our setting we have natural choices for such isomorphs given by whiskering: for instance, in the case of $\athree \otimes \atwo((0,0),(2,1))=A+A$ these are the sets $\{a(g,0):a \in A\}$ and $\{(h,1)a:a \in A\}$.\newline{}
We will use uniqueness of normal forms to show that $A$ is the trivial monoid.  Let $a \in A$ be a non-trivial element and consider the normal form decomposition $m^{2}_{h}a=b_{1} \ldots b_{n}$ which alternates between elements of the sets $(h,1)\{A\}$ and $\{A\}(g,0)$.  The maps $(i^{2}_{h},1)$ and $(1,i^{2}_{h})$ evaluate $b_{1}\ldots b_{n}$ to the product $\Pi_{j=2i}(a_{j})$ and $\Pi_{j=2i+1}(a_{j})$ respectively.  (Here $a_{j}$ is the $A$-component of $b_{j}$, and the empty product is identified with the unit of $A$.)  Since $(i^{2}_{h},1)m^{2}_{h}=1=(1,i^{2}_{h})m^{2}_{h}$ it follows that the normal form of $m^{2}_{h}a$ must have length at least $2$, as must $m^{2}_{v}a$ by an identical argument.\newline{}
Consider $m^{3}_{v} m^{2}_{h}a$.  Now $m^{3}_{v}=m^{2}_{v}+_{m^{1}_{v}}m^{2}_{v}$ maps words $(h,1)a$ and $a(g,0)$ of length $1$ to $(h,2)(m^{2}_{v}a)$ and $(m^{2}_{v}a)(g,0)$ respectively, whose normal forms are alternating words in $(h,2)A(0,g)$ and $(h,h)A$, and in $(2,h)A(g,0)$ and $A(g,g)$ respectively.  In particular each $m^{3}_{v} b_{i}$ is an alternating word of length $\geq 2$ in one pair, and $m^{3}_{v} b_{i+1}$ in the other pair.  Therefore the normal form of $m^{3}_{v}m^{2}_{h} a$ is simply obtained by juxtaposing the normal forms of the components $(m^{3}_{v} b_{1},\ldots ,m^{3}_{v} b_{n})$.  In particular the first two elements of this normal form must be either an alternating word in $(h,2)A(0,g)$ and $(h,h)A$, or in $(2,h)A(g,0)$ and $A(g,g)$: a \emph{vertically alternating} word.  By contrast, the first two elements of the normal form of $m^{3}_{h} m^{2}_{v}a$ alternate horizontally.  By uniqueness of normal forms we conclude that $A$ must be the trivial monoid.   \newline{}
Having shown that $\atwo \otimes \atwo$ is merely a category in the two possible cases, we observe that since the full inclusion $j:\Cat \to \SCat$ is closed under colimits and $\otimes$ cocontinuous in each variable, it follows that the pushouts $\athree \otimes \atwo$, $\atwo \otimes \athree$ and $\athree \otimes \athree$ must be categories too.
\epf
\begin{Theorem}\label{thm:TensorGray}
Let $\otimes:\GCat^{2} \to \GCat$ be a cocontinuous bifunctor with unit.  Then the underlying bifunctor
$$\cd{\Cat^{2} \ar[r]^{j^{2}} & \GCat^{2} \ar[r]^{\otimes} & \GCat \ar[r]^{(-)_{1}} & \Cat}$$
coincides with either the cartesian or funny tensor product.  Moreover, at categories $A$ and $B$, the tensor product $A \otimes B$ has only identity 2-cells.
\end{Theorem}
\pf
Since the forgetful functor $(-)_{2}:\GCat \to \SCat$ is cocontinuous the composite bifunctor $\otimes^{\prime}$:
$$\cd{\Cat^{2} \ar[r]^{j^{2}} & \GCat^{2} \ar[r]^{\otimes} & \GCat \ar[r]^{(-)_{2}} & \SCat}$$
is a cocontinuous bifunctor.  The full theorem amounts to the assertion that $\otimes^{\prime}$ is isomorphic to one of $j \circ \times,j \circ \star:\Cat^{2} \rightrightarrows \Cat \to \SCat$
which, by Proposition~\ref{prop:equiv}, is equally to show that the double cocategory $\arr \otimes^{\prime} \arr$ is isomorphic to $\arr \times \arr$ or $\arr \star \arr$.  By Remark~\ref{rk:HigherUnit} the unit for $\otimes$ is $1$ and it follows that $\arr \otimes^{\prime} \arr$ has left and top cocategories $\arr \otimes 1$ and $1 \otimes \arr$ isomorphic to $\arr$.  By Proposition~\ref{prop:sesquicat} $\arr \otimes^{\prime} \arr$ must be isomorphic to one of $\arr \star \arr$ and $\arr \times \arr$.
\epf

\subsection{No weak transformations.}
The informal argument of the introduction was against the existence of a monoidal biclosed structure on \GCat capturing weak transformations.  Here we make a precise statement. In Remark~\ref{rk:HigherUnit} we observed that any monoidal biclosed structure on $\GCat$ has unit $1$.  Accordingly the (left or right) internal hom $[A,B]$ must have Gray-functors as objects.  We take it as given that a weak transformation $\eta:F \to G$ of Gray-functors ought to involve components $\eta_{a}:Fa \to Ga$ and 2-cells $$\xy
(0,0)*+{Fa}="00"; (15,0)*+{Fb}="10";
(0,-15)*+{Ga}="01";(15,-15)*+{Gb}="11";
{\ar^{F\alpha} "00"; "10"};
{\ar_{\eta_{a}} "00"; "01"};
{\ar^{\eta_{b}} "10"; "11"};
{\ar_{G\alpha} "01"; "11"};
{\ar@{=>}^{\eta_{\alpha}}(11,-5)*+{};(5,-10)*+{}};
\endxy$$
(or in the opposite direction) to begin with.  Even if we restrict our attention to the simple case that $A$ is merely a category and $B$ a 2-category one expects that there should exist transformations whose components $\eta_{\alpha}$ are not identities.  This is, at least, the case for the pseudonatural transformations that arise in 2-category theory.  The following result shows that this is not possible.

\begin{Corollary}\label{cor:WeakTransformations}
Let $\otimes$ be a biclosed bifunctor with unit on $\GCat$, and let $[A,B]$ denote either internal hom.  Suppose that $A$ is a category and $B$ a 2-category.  Then the underlying category of $[A,B]$ is isomorphic to the category $[A,B_{1}]_{f}$ of functors and unnatural transformations, or to the category of functors and natural transformations $[A,B_{1}]$.
\end{Corollary}
\pf
%We will show that the underlying category $[A,B]_{1}$ of $[A,B]$ is either isomorphic to the category $[A,B_{1}]_{f}$ of functors and unnatural transformations or %isomorphic to the category $[A,B_{1}]$ of functors and natural transformations.\\
To begin, we observe that although the inclusion $j:\twocat \to \GCat$ does not have a right adjoint it does have a left adjoint $\Pi_{2}$.  This sends a Gray-category $A$ to the 2-category $\Pi_{2}(A)$ with the same underlying category as $A$ and whose 2-cells are equivalence classes of those in $A$: here, a parallel pair of 2-cells are identified if connected by a zig-zag of $3$-cells.  We will use the evident fact that if all 2-cells in $A$ are identities then $\Pi_{2}(A)$ is isomorphic to the underlying category $A_{1}$ of $A$.  We will also use the series of adjunctions
\begin{equation*}
\cd{\GCat \ar@/^1ex/[rr]^{\Pi_{2}} && \twocat \ar@/^1ex/[ll]^{j} \ar@/^1ex/[rr]^{(-)_{1}} && \Cat \ar@/^1ex/[ll]^{j}}
\end{equation*}
in which the left adjoints point rightwards, and in which the composite left adjoint is $(-)_{1}:\GCat \to \Cat$.\\
Now let $A$ be a category and $B$ a 2-category and $[A,B]$ the hom Gray-category characterised by the adjunction $- \otimes A \dashv [A,-]$.  We have isomorphisms:

\begin{gather*}
\Cat(C,([A,B])_{1}) \cong \GCat(C,[A,B]) \cong \GCat(C \otimes A,B) \cong \\
\twocat(\Pi_{2}(C \otimes A),B) \cong \twocat((C \otimes A)_{1},B)
\end{gather*}
in which the action of the inclusions $j$ are omitted.  The first three isomorphisms use the aforementioned adjunctions.  By Theorem~\ref{thm:TensorGray} all 2-cells in $C \otimes A$ are identities and this gives the fourth isomorphism.  Furthermore Theorem~\ref{thm:TensorGray} tells us that $(C \otimes A)_{1}$ is isomorphic to the funny tensor product of categories $C \star A$ or to the cartesian product $C \times A$.  In the first case we obtain further isomorphisms
\begin{gather*}
 \twocat((C \otimes A)_{1},B) \cong \twocat(C \star A,B) \cong \Cat(C \star A,B_{1}) \cong \Cat(C, [A,B_{1}]_{f}) \hspace{0.5cm} .
\end{gather*}
Applying the Yoneda lemma to the composite isomorphism $$\Cat(C,([A,B])_{1}) \cong  \Cat(C, [A,B_{1}]_{f})$$ we obtain $([A,B])_{1} \cong  [A,B_{1}]_{f}$ or, in the second case, an isomorphism $[A,B]_{1} \cong  [A,B_{1}]$. Finally, we note that an essentially identical argument works in the case that $[A,B]$ is the Gray-category corresponding to the adjunction $A \otimes - \dashv [A,-]$.
\epf

\subsection{No monoidal model structure.}
The cartesian product of categories and Gray tensor product of 2-categories are generally considered to be \emph{homotopically well behaved} tensor products.  In the former case we have that the cartesian product of equivalences of categories is again an equivalence; in the latter we have the analogous fact \cite{Lack2002A-quillen} but also further evidence, such as the fact that every tricategory is equivalent to a Gray-category \cite{Gordon1995Coherence}.  Still further evidence comes from the theory of Quillen model categories \cite{Quillen1967Homotopical}: both tensor products form part of \emph{monoidal model structures} \cite{Hovey1999Model}.  Recall that a monoidal model category consists of a category $\f C$ equipped with both the structure of a symmetric monoidal closed category and a model category; furthermore, these two structures are required to interact appropriately.  There is a unit condition and a condition called the pushout product axiom, a special case of which asserts that for each cofibrant object $A$ the adjunction $$A \otimes - \dashv [A,-]$$ is a Quillen adjunction.
% main condition is the pushout product axiom which asserts that for $f:A \to B$ and $g:C \to D$ the pushout product $$f \Box g:A \otimes D +_{A \otimes C} B \otimes C \to B \otimes D$$ is a cofibration if $f$ and $g$ are, and a trivial cofibration if, additionally, either $f$ or $g$ is.
When the unit is cofibrant, as it often is, this special case of the pushout product axiom already implies that the homotopy category $ho(\f C)$ is itself symmetric monoidal closed, and that the projection $\f C \to ho(\f C)$ is strong monoidal.\newline{}
In the case of \Cat the model structure in question has weak equivalences the equivalences of categories, whilst the fibrations and trivial fibrations are the isofibrations and surjective equivalences respectively.  In the case of $\twocat$ the weak equivalences are the biequivalences, whilst the fibrations are slightly more complicated to describe -- see \cite{Lack2002A-quillen}.  However each object is fibrant, and the trivial fibrations are those 2-functors which are both surjective on objects and locally trivial fibrations in \Cat.
 \newline{}
The model structure on \twocat is determined by that on \Cat in the following way.  Given a monoidal model category $\f V$ one can define a $\f V$-category $A$ to be \emph{fibrant} if each $A(a,b)$ is fibrant in $\f V$, and a $\f V$-functor $F:A \to B$ to be a \emph{trivial fibration} if it is surjective on objects and locally a trivial fibration in $\f V$.  Now a model structure is determined by its trivial fibrations and fibrant objects.  Consequently, when the above two classes do in fact determine a model structure on $\vcat$, the authors of \cite{Berger2013On-the-homotopy} referred to it as the \emph{canonical model structure} and described conditions under which it exists.  In particular the model structure on \twocat is canonical.\newline{}
Furthermore it was shown in \cite{Lack2011A-quillen} that \GCat admits the canonical model structure, as lifted from \twocat, and a natural question to ask is whether this forms part of a monoidal model structure.  %Our own interest in this question stemmed from the possibility that categories enriched in such a monoidal model structure might model all weak 4-categories.  
In 1.8(v) of \cite{Berger2013On-the-homotopy} the authors indicated that it was unknown whether such a tensor product exists.  The following result shows that this is too much to ask.
%\emph{Just wrote the above today and the last bit is a bit dull -- better off mentioning something about a monoidal model structure on \GCat providing a candidate for semistrict 4-categories.}
\begin{Corollary}\label{cor:model}
There exists no cocontinuous bifunctor $\otimes:\GCat^{2} \to \GCat$ with unit such that for each cofibrant $A$ either of the functors $A \otimes -$ or $- \otimes A$ is left Quillen.  In particular there exists no monoidal model structure on $\GCat$.
\end{Corollary}
\pf
By Corollary 9.4 of \cite{Lack2011A-quillen} a Gray-category is cofibrant just when its underlying sesquicategory is \emph{free on a computad}.  All we need are two immediate consequences of this: namely, that $\atwo$ is cofibrant and that each cofibrant Gray-category has \emph{underlying category free on a graph}.\newline{}
Now suppose such a bifunctor $\otimes$ does exist.  Since $\atwo$ is cofibrant one of $\atwo \otimes -$ and $- \otimes \atwo$ must be left Quillen.  Therefore $\atwo \otimes \atwo$ is cofibrant.  In particular the underlying category of $\atwo \otimes \atwo$ must be free on a graph.  Now $(\atwo \otimes \atwo)_{1} = \atwo \star \atwo$ or $\atwo \times \atwo$ by Theorem~\ref{thm:TensorGray} and only the first of these is free on a graph.  Consequently $(\atwo \otimes \atwo)_{1} = \atwo \star \atwo$.\newline{}
Let $I$ denote the free isomorphism  $0 \cong 1$ and consider the map $j:\atwo \to I$ sending $0 \to 1$ to $0 \cong 1$.  By Theorem~\ref{thm:TensorGray} $j \otimes j:\atwo \otimes \atwo \to I \otimes I$ coincides with the funny tensor product $j \star j:\atwo \star \atwo \to I \star I$ on underlying categories, and both 2-categories have only identity 2-cells.  In particular $j \otimes j$ sends the non-commuting square
$$\xy
(0,0)*+{(0,0)}="00"; (15,0)*+{(1,0)}="10";(0,-15)*+{(0,1)}="01";(15,-15)*+{(1,1)}="11";
{\ar_{(0,f)} "00"; "01"};
{\ar^{(f,0)} "00"; "10"};
{\ar^{(1,f)} "10"; "11"};
{\ar_{(0,f)} "01"; "11"};
\endxy$$
to a \emph{still non-commuting} square of isomorphisms.\newline{}
As described in Section 3 of \cite{Lack2011A-quillen} there is a \emph{universal adjoint biequivalence} $E$ in \GCat.  It is a cofibrant Gray-category with two objects and both maps $1 \rightrightarrows E$ are trivial cofibrations.  In particular $E$ has underlying category free on the graph $(f:0 \leftrightarrows 1:u)$ so that we obtain a factorisation of $j$:
\[
\atwo \stackrel{j_{1}}{\longrightarrow} E \stackrel{j_{2}}{\longrightarrow} I \rlap{ ,}
\]
and hence of $j \otimes j$:
\[
\atwo \otimes \atwo \stackrel{l}{\longrightarrow} E \otimes E \stackrel{k}{\longrightarrow} I \otimes I \rlap{ ,}
\]
where we have written $j_{1} \otimes j_{1} = l, j_{2} \otimes j_{2} = k$.  The two distinct paths $r,s:(0,0) \rightrightarrows (1,1)$ of $\atwo \otimes \atwo$ give rise to parallel 1-cells $l(r),l(s):l(0,0) \rightrightarrows l(1,1)$.  A 2-cell $\alpha:l(r) \Rightarrow l(s)$ would yield a 2-cell $k(\alpha):j \otimes j(r) \Rightarrow j\otimes j(s)$ but no such 2-cell can exist by the above.  In particular there exists no 2-cell $l(r) \Rightarrow l(s)$ in $E$.\newline{}
However since $1 \to E$ is a trivial cofibration in \GCat we must have that $1 \cong 1 \otimes 1 \to 1 \otimes E \to E \otimes E$ is a trivial cofibration.  By 2 from 3 the map $E \otimes E \to 1$ must be a weak equivalence and indeed, since all Gray-categories are fibrant, a trivial fibration.  But this implies that there exists a 2-cell between any parallel pair of 1-cells in $E \otimes E$ and we have just shown that this is not the case.
\epf

\subsection{Dolan's objection to Street's files.}
In a 1996 email \cite{Dolan1996Email} to Ross Street, James Dolan described an objection to Street's then conjectural notion of semistrict n-category -- called an \emph{n-file} \cite{Street1995Notes}.  Dolan's objection, brought to our attention during the editorial process, is closely related to the arguments of the present paper and for that reason we discuss it here.\\
The category $n\textnormal{-Cat}$ of strict $n$-categories is cartesian closed.  Taking the union of the sequence $\textnormal{(n+1)-Cat} = (n\textnormal{-Cat},\times)\textnormal{-Cat}$ gives the cartesian closed category $\f V_{1}=(\OCat,\times)$.\begin{footnote}{As in \cite{Street1995Notes,Street2004Categorical-and} an object of $\OCat$ is an $n$-category for some $n$.  We note that $\OCat$ is often used to refer to $\omega$-categories with non-trivial cells in each finite dimension.}\end{footnote}  $\OCat$ supports several other monoidal biclosed structures, including a higher dimensional version of the (lax) Gray tensor product.  The following construction of this tensor product follows \cite{Street2004Categorical-and} and we refer to that paper for further details and references.  The free $\omega$-categories $\{O(I^{n}):n \in \mathbf N\}$ on the \emph{parity $n$-cubes $I^{n}$} form a dense full subcategory $Q$ of $\OCat$.  The subcategory is monoidal, with tensor product $O(I^{n}) \otimes O(I^{m}) = O(I^{n+m})$ induced from the tensor product of cubes.  Day's technique \cite{Day1970On-closed, Day1972A-reflection} for extending monoidal structures along a dense functor can, with work, be employed to establish the monoidal biclosed structure $\f V_{2}=(\OCat,\otimes)$: here called the Gray-tensor product of $\omega$-categories.  Gray-categories, identifiable as 3-dimensional $\f V_{2}$-categories with invertible coherence constraints, were called \emph{$3$-files} in \cite{Street1995Notes}.\\
Corresponding to the fact that the cartesian product is the terminal semicartesian tensor product, so $\f V_{1}$-categories, just $\omega$-categories, can naturally be viewed as $\f V_{2}$-categories.  Accordingly we have a composite functor
\begin{equation*}
\cd{Q \ar[r] & \OCat \simeq \f V_{1}\textnormal{-Cat} \ar[r] & \f V_{2}\textnormal{-Cat}}
\end{equation*}
and the composite is fully faithful.  It was conjectured in \cite{Street1995Notes} that this functor satisfied Day's conditions whereby the monoidal structure on $Q$ could be extended to a biclosed one $\otimes_{2}$ on $\f V_{2}\textnormal{-Cat}$; then one would take $\f V_{3}=(\f V_{2}\textnormal{-Cat},\otimes_{2})$, define 4-files as certain 4-dimensional categories enriched in $\f V_{3}$, and iterate to higher dimensions.\\
Dolan's objection\begin{footnote}{In \cite{Dolan1996Email} Dolan writes ``The category $Q$ of Gray tensor powers of the free strict infinity category $I$ on one 1-arrow is a full subcategory of $V_{2}$, and also of $V_{3}$.  It's also a monoidal subcategory of $V_{2}$, and you want it to be isomorphically a monoidal subcategory of $V_{3}$ as well.  But then since $I$ is a co-category object in $V_{3}$ (just as it is in $V_{2}$), $I*I$ would have to be a double co-category object in $V_{3}$, with moreover several (4, at least) co-category homomorphisms from $I$.  Exhaustive analysis of $I*I$ shows that there is only one plausible candidate for the horizontal co-category structure on it, and only one plausible candidate for the vertical co-category structure on it.  But these two co-category structures do not strictly commute; there is the contradiction."}\end{footnote}, expressed in the terms of the present paper, was as follows.  The co-graph $O(I^{0})\rightrightarrows O(I^{1}) = d,c:1 \rightrightarrows \atwo$ underlies the standard cocategory $\arr$ in $\Cat$, equally a cocategory in $\OCat$ or in $\f V_{2}\textnormal{-Cat}$.  Accordingly a biclosed tensor product $\otimes_{2}$ would give a double cocategory $\arr \otimes_{2} \arr$ with left and top cocategories $\arr \otimes_{2} 1$ and $1 \otimes_{2} \arr$ isomorphic to $\arr$.  If $\otimes_{2}$ were obtained by Day's technique of extension then the composite inclusion $Q \to \f V_{2}\textnormal{-Cat}$ would be strong monoidal; hence $(\arr \otimes_{2} \arr)_{1,1}=O(I^{1}) \otimes_{2} O(I^{1}) \cong O(I^{2})$.  But $O(I^{2})$ is just 
$$\xy
(0,-7.5)*+{(0,0)}="00"; (15,-7.5)*+{(1,0)}="10";
(0,-22.5)*+{(0,1)}="01";(15,-22.5)*+{(1,1)}="11";(23,-24)*+{.}="";
{\ar^{} "00"; "10"};
{\ar_{} "00"; "01"};
{\ar^{} "10"; "11"};
{\ar_{} "01"; "11"};
{\ar@{=>}^{}(11,-12.5)*+{};(5,-17.5)*+{}};
\endxy$$
Dolan observed that there is only one possible double cocategory structure on $O(I^{2})$ (with $\arr$ as its left and top cocategories) and that this does not form a double cocategory: as in the $\GCat$ setting of Section 4, the middle four interchange axiom fails to hold.  Consequently there cannot exist a biclosed tensor product $\otimes_{2}$ obtained by extension along $Q \to \f V_{2}\textnormal{-Cat}$.\\
The key observation above is that the lax square double cocategory in $\twocat$ (or $\OCat$) does not remain one on passing to a world like $\GCat$ (or $\f V_{2}\textnormal{-Cat}$) in which the middle 4-interchange doesn't hold.  This is the same observation at the heart of our key result Proposition~\ref{prop:sesquicat} -- whose additional value is that it examines not just the lax square but all possible double cocategories.


\begin{thebibliography}{10}

\bibitem{Adamek1994Locally}
{\sc Ad{\'a}mek, J., and Rosick{\'y}, J.}
\newblock Locally presentable and accessible categories.
\newblock {\em London Math. Soc. Lecture Notes 189}
\newblock xiv+316 pp., Cambridge University Press (1994)

\bibitem{Berger2013On-the-homotopy}
{\sc Berger, C., and Moerdijk, I.}
\newblock On the homotopy theory of enriched categories.
\newblock {\em Quarterly Journal of Mathematics 64}, 3 (2013), 805--846.

\bibitem{Bergman1998An-invitation}
{\sc Bergman, G.}
\newblock {\em An invitation to general algebra and universal constructions}.
\newblock Henry Helson, Berkeley, CA, 1998.

\bibitem{Crans1999A-tensor}
{\sc Crans, S.}
\newblock A tensor product for {G}ray-categories.
\newblock {\em Theory and Applications of Categories 5}, 2 (1999), 12--69.

\bibitem{Day1970On-closed}
{\sc Day, B.}
\newblock On closed categories of functors.
\newblock In {\em Midwest Category Seminar Reports IV}, vol.~137 of {\em
  Lecture Notes in Mathematics}. (Springer-Verlag, Berlin 1970) pp.~1-38.
 
 \bibitem{Day1972A-reflection}
{ \sc Day, B.}
 \newblock A reflection theorem for closed categories.
 \newblock{\em Journal of Pure and Applied Algebra 2}, 1 (1972), 1-11.

\bibitem{Day1997Monoidal}
{\sc Day, B., and Street, R.}
\newblock Monoidal bicategories and {H}opf algebroids.
\newblock {\em Advances in Mathematics 129}, 1 (1997), 99--157.

\bibitem{Dolan1996Email}
{\sc Dolan, James.}
\newblock{n-files.}
\newblock{\em Email to Ross Street} (March 1996).

\bibitem{Ehresmann1963Categories}
{\sc Ehresmann, C.}
\newblock{Cat\'{e}gories structur{\'e}es. III.  Quintettes et applications covariantes}.
\newblock{ In \em{Topol. et. G{\'e}om. Diff.} (S\'em. C. Ehresmann), Vol. 5, page 21}.
\newblock{Institut H. Poincar\'{e} Paris, 1963.}

\bibitem{Foltz1980Algebraic}
{\sc Foltz, F., Lair, C., and Kelly, G.~M.}
\newblock Algebraic categories with few monoidal biclosed structures or none.
\newblock {\em Journal of Pure and Applied Algebra 17}, 2 (1980), 171--177.

\bibitem{Gordon1995Coherence}
{\sc Gordon, R., Power, J., and Street, R.}
\newblock Coherence for tricategories.
\newblock {\em Memoirs of the American Mathematical Society 117\/} (1995).

\bibitem{Gray1974Formal}
{\sc Gray, J.}
\newblock Formal category theory: adjointness for 2-categories.
\newblock {vol.~391 of Lecture Notes in Mathematics  Springer, 1974.}

\bibitem{Hovey1999Model}
{\sc Hovey, M.}
\newblock {\em Model categories}, vol.~63 of {\em Mathematical Surveys and
  Monographs}.
\newblock American Mathematical Society, 1999.

\bibitem{Kelly1974Doctrinal}
{\sc Kelly, G.~M.}
\newblock Doctrinal adjunction.
\newblock In {\em Category Seminar (Sydney, 1972/1973)}, vol.~420 of {\em
  Lecture Notes in Mathematics}. Springer, 1974, pp.~257--280.

\bibitem{Lack2002A-quillen}
{\sc Lack, S.}
\newblock A {Q}uillen model structure for 2-categories.
\newblock {\em K-Theory 26}, 2 (2002), 171--205.

\bibitem{Lack2011A-quillen}
{\sc Lack, S.}
\newblock A {Q}uillen model structure for {G}ray-categories.
\newblock {\em Journal of K-Theory 8}, 2 (2011), 183--221.

\bibitem{Lumsdaine2011A-small}
{\sc Lumsdaine, P.~L.}
\newblock A small observation on co-categories.
\newblock {\em Theory and Applications of Categories 25}, 9 (2011), 247--250.

\bibitem{Quillen1967Homotopical}
{\sc Quillen, D.~G.}
\newblock {\em Homotopical algebra}, vol.~43 of {\em Lecture Notes in
  Mathematics}.
\newblock Springer, 1967.

\bibitem{Pare1974Colimits}
{\sc Pare, R.}
\newblock Colimits in Topoi.
\newblock{ \em Bull. Amer. Math. Soc. 80}, (1974), 556--561.

\bibitem{Street1995Notes}
{\sc Street, Ross.}
{\newblock Descent theory. Notes from Oberwolfach lectures of 17-23 December 1995.}
{{Available at \tt http://maths.mq.edu.au/~street/Publications.htm}.}

\bibitem{Street1996Categorical}
{\sc Street, Ross.}
{\newblock Categorical Structures,}
{\newblock Handbook of algebra {V}ol. 1}.
\newblock North-Holland, 1996, pp.~529--277.

\bibitem{Street2004Categorical-and}
{\sc Street, Ross.}
{\newblock Categorical and combinatorial aspects of descent theory,}
\newblock Appl. Categ. Structures 12, (2004), 5-6, 537-576.


\bibitem{nlab2009} {\tt http://ncatlab.org/nlab/revision/generalized+Gray+tensor+product/1} (2009).


\end{thebibliography}
\end{document}